\documentclass[12 pt]{article}
\usepackage{amssymb,amsmath,amsfonts,amsthm}
\usepackage{xcolor}
\newcommand\R{{\mathbb {R}}}
\renewcommand\P{{\mathbb {P}}}

\newcommand\N{{\mathbb {N}}}

\newcommand\E{{\mathbb {E}}}

\newtheorem {Theorem}{Theorem}[section]
\newtheorem {Lemma}[Theorem]{Lemma}
\newtheorem {Proposition}[Theorem]{Proposition}

\theoremstyle{definition}
\newtheorem{Definition}{Definition}[section]

\newcommand\beq{\begin{equation}}
\newcommand\eeq{\end{equation}}

\begin{document}

\title{Limit theorems  for iid products of positive matrices}

\author{C. Cuny\footnote{Christophe Cuny, Univ Brest, UMR 6205 CNRS 6205, LMBA, 6 avenue Victor Le Gorgeu, 29238 Brest}, J. Dedecker\footnote{J\'er\^ome Dedecker, Universit\'e de Paris, CNRS, MAP5, UMR 8145,
45 rue des  Saints-P\`eres,
F-75006 Paris, France.}
and 
 F. Merlev\`ede \footnote{Florence Merlev\`ede, LAMA,  Univ Gustave Eiffel, Univ Paris Est Cr\'eteil, UMR 8050 CNRS,  \  F-77454 Marne-La-Vall\'ee, France.}}

\maketitle

\begin{abstract} 
We study stochastic properties of the norm cocycle associated with iid products of positive matrices. We obtain the almost sure invariance principle (ASIP) with rate $o(n^{1/p})$ under the optimal condition  of a moment or order $p>2$ and the Berry-Esseen theorem with rate $O(1/\sqrt n)$ under the optimal condition of a moment of order 3.  The results are also valid for the matrix norm. For the matrix coefficients, we also have the ASIP but we obtain only partial results for the Berry-Esseen theorem. The proofs make use of coupling coefficients that surprisingly decay exponentially fast to 0 while there is only a polynomial decay in the case of invertible matrices.  
All the results are actually valid in the context of iid products of matrices leaving invariant a suitable cone.
\end{abstract}

{\it AMS 2020  subject classifications}: 60F05, 60B15, 60G50. 

{\it Key words and phrases}. Random walk;  Cocycle;  Berry-Esseen theorem, almost sure invariance principle, Hilbert metric.

\section{Introduction}

In  a series of paper \cite{CDJ}, \cite{CDM-deviation}, \cite{CDM}, \cite{CDMP} and \cite{CDMP-bis} we studied the stochastic properties of the norm cocycle associated with the 
left random walk on $GL_d(\R)$ under optimal or close to optimal moment conditions. The moment conditions are optimal in case of the central limit theorem (CLT) and the ASIP with rate and close to optimal in the case of the 
Berry-Esseen theorem. We also obtained results for the matrix norm, the matrix coefficients and the spectral radius.

\medskip

A key ingredient in the proofs is the use of some coupling 
coefficients introduced in \cite{CDJ}, see Section \ref{section-coupling} for the definition.

\medskip

It turns out that it is also possible to control similar coefficients in the context of the left random walk on 
the semi-group of matrices of size $d\ge 2$, with non-negative entries (that we call positive matrices in the sequel). Actually, one can even prove the exponential convergence to 0 of those coefficients under polynomial moment conditions, see Proposition \ref{delta}. As a  consequence, we obtain 
Berry-Essen's theorem with rate $O(1/\sqrt n)$ under the optimal condition of a moment of order 3. We also obtain optimal intermediary rates under moments of order 
$p\in (2,3)$. Finally, we also obtain optimal rates in the ASIP.

\medskip

Let us mention that the study of iid products of positive matrices benefited from a lot of works. Let us cite, among others, Hennion \cite{Hennion}, Buraczewski et al. \cite{BDGM}, Buraczewski and Mentmeier \cite{BM} or Grama, Liu and Xiao
\cite{GLX-first}.

\medskip

Hennion obtains the strong law of large numbers and the CLT 
under optimal moment conditions in the more general situations of product of dependent positive random matrices, under mixing conditions. All the other above mentionned papers assume 
exponential moment which allows to use in a natural way 
the Guivarc'h-Nagev method, which is based on perturbation of operators. 

\medskip


It has been observed in the preprint \cite{GLX}, that the Guivarc'h-Nagaev method applies under polynomial moment 
conditions. In particular, they obtain the Berry-Esseen theorem with rate $O(1/\sqrt n)$ under a moment of order 3 plus some extra technical condition, see their condition (A2).

\medskip

In Section 2, we introduce some notations and definitions and we also recall several key properties in the study of positive matrices. 

\medskip

In section 3, we establish the existence of a unique invariant probability   and we estimate our coupling coefficents. 

\medskip

In section 4, we recall the strong law of large numbers of Hennion and provide some complementary results.

\medskip

In section 5, we recall the CLT and provide several identification of the asymptotic variance $s^2$. Moreover, we show that the known aperiodicity condition (see Definition \ref{aperiodic}) is sufficient for $s^2>0$, under a moment of order 2.

\medskip

In section 6, we obtain the ASIP for the norm cocycle, the matrix norm, the spectral radius and the matrix coefficients under optimal polynomial moment condition. We also consider the case of exponential moments, but we have a slight loss compare to the known result in the iid case (which corresponds to $d=1$ in our setting).

\medskip

 In section 7, we obtain  the Berry-Esseen theorem for all the above mentionned quantities. The obtained rates are optimal (in terms of moment conditions) in the case of the norm cocyle and the matrix norm, but we have a  loss in the case of the spectral radius and the matrix coefficients.

\medskip

In section 8 we study the regularity of the invariant measure and in section 9, we provide some deviation inequalities for the norm cocycle and the matrix coefficients.

\medskip

In section 10, we explain how to generalize our results to matrices leaving invariant a suitable cone (notice that the positive matrices of size $d$ may be seen as the matrices leaving invariant the cone 
$(\R^+)^d$. 

\medskip

Finally, in section 11, we provide technical results relevant to the previous section. 

\medskip

In all the paper we denote $\N:=\{1,2, \ldots\}$.

\section{Norm cocycle and matrix norm}

\setcounter{equation}{0}

Let $d\ge 2$ be an integer. Let $G$ be the semi-group of 
$d$-dimensional positive allowable matrices:  by positive, we mean that all entries are greater than or equal to 0, by allowable, we mean that  any lign and any column admits a strictly positive element.

\medskip

We endow $\R^d$ with the $\ell^1$ norm and $G$ with the corresponding operator norm. We denote both norms by $\|\cdot \|$. Recall that $\|g\|=\sup_{\|x\|=1}
\|gx\|$.

\medskip

We put on $G$ the topology inherited from (the distance associated with) the norm. Then, $G$ becomes a locally compact space.

\medskip

 Let $G^{+}$ be the sub-semi-group of $G$ whose entries are all strictly positive. Actually, $G^{+}$ is the interior of $G$.

\medskip

Define 
\begin{gather}\label{S+-1}
S^+:=\{x=(x_1, \ldots, x_d)\in \R^d\, :\, \|x\|=1\, \mbox{and } x_i\ge
 0, \, \forall i\in \{1,\ldots, d\}\, \}\, ,\\
 \label{S+-2}S^{++}:=\{x=(x_1, \ldots, x_d)\in \R^d\, :\, \|x\|=1\, \mbox{and } x_i>
 0, \, \forall i\in \{1,\ldots, d\}\, \}\, .
\end{gather}

\medskip

Notice that for $g\in G$, we actually have $\|g\|=\sup_{x\in S^+}\|gx\|$ and that, if $g=(g_{ij})_{1\le i,\, j\le d}$, 
\begin{equation}\label{norm}
 \|g\|=\max_{1\le j\le d}\sum_{i=1}^dg_{ij}\, .
\end{equation} 

\medskip

For every $g\in G$, set $v(g) =\inf_{x\in S^+} \|gx\|$. If $g=(g_{ij})_{1\le i,\, j\le d}$, we have
\begin{equation}\label{v-def}
v(g)=\min_{1\le j\le d}\sum_{i=1}^dg_{ij}\, .
\end{equation}   

By definition of $G$, $v(g)>0$ for every $g\in G$.

\medskip

We then define $N(g):=\max(\|g\|, 1/v(g))$ and 
$L(g)=\frac{\|g\|}{v(g)}$. Notice that $N(g)^2\ge L(g)\ge 1$ for every $g\in G$.

\medskip

We endow $S^+$ with the following metric (see Proposition \ref{complete-metric} for a proof that it is indead a metric). For every $x,\, y\in S^+$, 
$$
d(x,y) =\varphi(m(x,y)m(y,x))\, ,
$$
where 
\begin{equation}\label{phi}
\varphi(s)= \frac{1-s}{1+s}\qquad \forall s\in [0,1]\, ,
\end{equation}
and 
$$
m(u,v)=\inf\left\{\frac{u_i}{v_i}\, :\, i\in \{1,\ldots ,d\},\, v_i>0\, .\right\}
$$

Notice that the diameter of $S^+$ is 1 and that $d(x,y)=1$ if and only if there exists $i_0\in \{1, \ldots , d\}$ such that $x_{i_0}=0$ and $y_{i_0}>0$ or $x_{i_0}>0$ and $y_{i_0}=0$.


\medskip

Using that for $u,v\in S^+$, $\max_{1\le i\le d}u_i\le 1$ and $\max_{1\le i\le d}v_i\ge 1/d$, we see that $m(u,v)\le d$.

\medskip

The semi-group $G$ is acting  on $S^+$ as follows. 
$$
g\cdot x =\frac{gx}{\|gx\|} \qquad \forall (g,x)\in G\times 
S^+\, .
$$

We then define a cocyle by setting $\sigma(g,x)=
\log (\|gx\|)$ for every $(g,x)\in G\times S^+$. The cocycle property reads 
\begin{equation}\label{cocycle-property}
\sigma(gg',x)=\sigma(g, g'\cdot x)+\sigma(g',x)\, .
\end{equation}

\medskip

Following Hennion \cite[Lemma 10.6]{Hennion}, for every 
$g\in G$ we define $c(g):=\sup_{x,y\in  S^+}d(gx,gy)$.

\medskip

Let us recall some properties that one may find in 
Hennion \cite{Hennion}, see his Lemmas 5.2, 5.3 and 10.6 and his Proposition 3.1.

\begin{Proposition}\label{prop-hennion}
For every $(g,g',x,y)\in G^2\times(S^+)^2$ we have
\begin{itemize}
\item [$(i)$] $|\sigma(g,x)|\le \log N(g)$;
\item [$(ii)$] $\|x-y\|\le 2d(x,y)$; 
\item [$(iii)$] $|\sigma(g,x)-\sigma(g,y)|\le 2L(g)d(x,y)$;
\item [$(iv)$] $|\sigma(g,x)-\sigma(g,y)|\le 2\ln \big(1/
(1-d(x,y))\big)$ ;
\item [$(v)$] $c(gg')\le c(g)c(g')$;
\item [$(vi)$] $c(g)\le 1$ and $c(g)< 1$ 
iff $g\in G^+$;
\item [$(vii)$] $d(g\cdot x,g\cdot y)\le c(g) d(x,y)$.
\end{itemize}
\end{Proposition}

Let us also mention a closed-form expression for $c(g)$ obtained in Lemma 10.7 of \cite{Hennion}. For every $g=(g_{ij})_{1\le i,\, j\le d}$ we have 
\begin{equation}\label{closed}
c(g)=\max_{1\le i,\, j,\, k,\, \ell\le d}
\frac{|g_{ij}g_{k\ell}-g_{i\ell}g_{kj}|}{g_{ij}g_{k\ell}+g_{i\ell}g_{kj}}\, .
\end{equation}

Notice that $(g,x)\to g x$ is  continuous on $G\times S^+$ (for the distance on $G$ induced by the operator norm and the distance on $S^+$ induced by $\|\cdot \|$) and does not vanish. Hence, it follows from item $(ii)$ that $(g,x)\to g\cdot x$ is continuous on $G\times S^+$ (for the distance on $G$ induced by the operator norm and the distance $d$ on 
$S^+$).

\medskip

Let us give some more properties that will be useful in the sequel. Set $e=\{1/d, \ldots , 1/d\}\in S^+$. For $g\in G$, we denote by $g^t$ 
the adjoint matrix of $g$.

\begin{Lemma}\label{properties}
For every $(g,x,y)\in G\times (S^+)^2$,
\begin{itemize}
\item [$(i)$] $|\sigma(g,x)-\sigma(g,y)|\le \log L(g)$;
\item [$(ii)$] $\|g e\|\le \|g\|\le d\|ge\|$; 
\item [$(iii)$] $\|g\|\le d\|g^t\|$;
\item [$(iv)$] $|\sigma(g,x)-\sigma(g,y)|\le 2(2 +\log L(g)) 
d(x,y)$.
\end{itemize}
\end{Lemma}

\noindent {\bf Remark.} The inequality in item $(iv)$ of Lemma \ref{properties} is much better that the one in item $(iii)$ of Proposition \ref{prop-hennion}. 

\medskip

\noindent {\bf Proof.} Items $(i)$ and $(ii)$ are obvious.
Item $(iii)$ is an easy consequence of \eqref{norm}. Let us prove item 
$(iv)$. Let $x,\, y\in S^+$. Assume that $d(x,y)\le 1/2$. 
Notice that for every $t\in [0,1/2]$, $\ln (1/(1-t))\le 2t$. Hence, using item $(iv)$ of Proposition 
\ref{prop-hennion}, we see that $|\sigma(g,x)-\sigma(g,y)|\le 4 d(x,y)$. If $2d(x,y)\ge 1$, then the desired conclusion follows from item $(i)$ of Lemma \ref{properties}. \hfill $\square$

\medskip

\begin{Proposition}\label{complete}
$(S^+,d)$ is complete and $S^{++}$ is closed. 
\end{Proposition}

\noindent {\bf Remark.} A Hint of proof of completeness is given after Theorem 4.1 of Bushell \cite{Bushell}, for Hilbert's metric given by $d_H(x,y)=-\ln (m(x,y)m(y,x))$. See Proposition \ref{complete-metric} for a proof in a more general situation.

\medskip



 
 


Let us state some of the assumptions used throughout the paper.

\begin{Definition}
Let $\mu$ be a Borel probability on $G$ and $p\ge 1$.  We say that $\mu$ 
admits a moment of order $p$ if 
$$
\int_G (\log (N(g)))^pd\mu (g) <\infty\, .
$$
We say that $\mu$ almost admits a moment of order $p$ if
$$
\int_G (\log (L(g)))^pd\mu (g) <\infty\, .
$$
\end{Definition}

\noindent {\bf Remark.} Clearly, since $L(g)\le N(g)^2$, if $\mu$ admits a moment of order $p\ge 1$, it almost admits a moment of order $p\ge 1$, but the converse is not true in general, see the example in Section 6.  Assume now that $\mu$ almost admits a moment of order $p\ge 1$. Then, $\mu $ admits a moment of order $p$ iff $\int_G |\log \|g\| |^p d\mu(g)<\infty$ iff
$\int_G|\log v(g)|^pd\mu(g) <\infty$. 

\medskip

Similarly, we say that $\mu$ admits or almost admits an 
exponential moment of order $\gamma>0$, if there exists $\delta>0$ such that, respectively,
$$
\int_G {\rm e}^{\delta N(g)^\gamma} d\mu(g)<\infty\, ,
$$ 
or 
$$
\int_G {\rm e}^{\delta L(g)^\gamma} d\mu(g)<\infty\, .
$$

\begin{Definition}\label{contracting}
We say that $\mu$ is \emph{strictly contracting} if 
there exists $r\in \N$, such that $\mu^{*r}(G^+)>0$.
\end{Definition}

Equivalently, the closed semi-group $\Gamma_\mu$ generated by the support of $\mu$ has non empty intersection with 
$G^+$.

\medskip

\section{Invariant measure and coupling coefficients}\label{section-coupling}

\setcounter{equation}{0}

Recall that a Borel (with respect to $d$) probability $\nu$ on $S^+$ is said to 
be $\mu$-invariant if for every Borel non negative function $\varphi$ on $S^+$, $\int_{G\times S^+}\varphi(g\cdot x)
d\mu(g)d\nu(x)=\int_{S^+}\varphi(x)d\nu(x)$. It is 
well known and easy to prove (recall that 
$(g,x)\to g\cdot x)$ is continuous on $G\times S^+$) that the support of a $\mu$-invariant measure is $\Gamma_\mu$-invariant, i.e. satisfies $\Gamma_\mu\cdot {\rm supp}
\, \nu\subset {\rm supp}\, \nu$ .

\medskip

We will see that when $\mu$ is strictly contracting, it admits a unique $\mu$-invariant probability on $S^+$. 
We need some further notation to identify its support. 

\medskip

Let $g\in G^+$. By the Perron-Frobenius theorem (see Theorem 1.1.1 of \cite{LN}), there exists a unique $x\in S^{++}$ 
such that $gx=\kappa(g) x$, where $\kappa(g)$ is the spectral radius of $g$. We denote that vector by $v_g$. 
Then, clearly, we have 
\begin{equation}\label{v-radius}
\kappa (g)\ge v(g)\qquad \forall g\in G\, .
\end{equation}


\medskip

\medskip

Following \cite{BDGM} (see (2.4) there) we define
$$
\Lambda_\mu=\overline{\{v_g\, :\, g\in \Gamma_\mu\cap G^+\}}\, ,
$$
where the closure is taken with respect to $d$. By Proposition \ref{complete}, $\Lambda_\mu\subset S^{++}$.

It follows from Lemma 4.2 of \cite{BDGM} that 
$\Lambda_\mu$ is $\Gamma_\mu$-invariant (i.e. 
$\Gamma_\mu\cdot \Lambda_\mu\subset \Lambda_\mu$).


\medskip

 The existence and uniqueness in the next proposition follow from Theorem 2.1 of \cite{HH}. We provide a slightly different proof and identify the support of the invariant measure.

\begin{Proposition}
Assume that $\mu$ is strictly contracting. Then, 
there exists 
a unique $\mu$-invariant  probability $\nu$ on $S^+$. Moreover ${\rm supp}\, \nu= \Lambda_\mu$.
\end{Proposition}


\noindent {\bf Proof.} Let $(Y_n)_{n\in \N}$ be iid random variables taking values in $G$, with law $\mu$. Let $r\in \N$ be as in Definition 
\ref{contracting}. For every  $n\in \N$, set $B_n:= Y_1\cdots Y_n$. Let $m:=[n/r]$. Notice that, by item $(v)$ of Proposition \ref{prop-hennion}, 
$c(B_n)\le \prod_{k=0}^{m-1}c(Y_{kr+1}\cdots Y_{(k+1)r})$. 
By the strong law of large numbers and the fact that $\mu$ is strictly contracting, using item $(vi)$ of Proposition \ref{prop-hennion},
$$\frac1m \sum_{k=0}^{m-1}\log c(Y_{kr+1}\cdots Y_{(k+1)r})\underset{m\to +\infty}\longrightarrow \E(\log c(Y_{1}\cdots Y_{r}))<0\qquad  \mbox{$\P$-a.s.}
$$
 Hence, $c(B_n)=O(\delta^m)$ 
almost surely for some $0<\delta <1$. In particular, $c(B_n)<1$ for $n$ large enough, so that, by item $(vi)$ of Proposition \ref{prop-hennion}, $B_n\in G^+$ and $B_n\cdot x\in S^{++}$ for every $x\in S^+$.

\medskip

Let $x\in S^+$. By item $(vii)$ of Proposition \ref{prop-hennion}, there exists a non negative random variable $K$, independent of $x$, such that for every $n\in \N$,
$$
d(B_n\cdot x,B_{n+1}\cdot x)\le c(B_n)\le K \delta^m\, .
$$
Hence $(B_n\cdot x)_{n\in \N}$ is Cauchy, taking values in $S^{++}$ for $n$ large enough, hence converges 
to some random variable $Z$ whose law is $\mu$-invariant. 
By item $(vii)$ of Proposition \ref{prop-hennion}, $d(B_n\cdot x,B_n\cdot y)\le c(B_n)$ and we see that 
$(B_n\cdot y)_{n\in \N}$ converges to $Z$ for every 
$y\in S^+$.

\medskip

Let $\nu$ be a $\mu$-invariant probability on $S^+$. 
Then, for every $m\in \N$, and every continuous bounded 
$\varphi$ on $S^+$, we have 
$$
\int_{S^+}\varphi d\nu= \int_{S^+}\E \big(\varphi (B_{m}\cdot x)\big) d\nu(x)\underset{m\to +\infty}\longrightarrow 
\E(\varphi(Z))\, ,
$$
which proves uniqueness.

\medskip

The fact that ${\rm supp}\, \nu\supset \Lambda_\mu$ 
follows from the fact that ${\rm supp} \, \nu $
is $\Gamma_\mu$-invariant and from Lemma 4.2 of 
\cite{BDGM}. To prove the converse inclusion, just notice that, since $\Gamma_\mu \cdot \Lambda_\mu\subset \Lambda_\mu $, 
for every $x\in \Lambda_\mu$, $B_n\cdot x\in \Lambda_\mu$ 
almost surely. Hence $Z\in \Lambda_\mu$ almost surely 
and $\nu(\Lambda_\mu)=1$ which implies the desired result. 
\hfill $\square$

\medskip

Let $(Y_n)_{n\in \N}$ be iid random variables taking values in $G$, with law $\mu$. For every $n\in \N$, set $A_n:=Y_n\cdots Y_1$.

\medskip

For every $p\ge 1$ and every $n\in \N$ define 
$$
\delta_{p,\infty}(n):=\sup_{x,y\in S^+}
\E\big(|\sigma(Y_n, A_{n-1}\cdot x)-\sigma(Y_n,A_{n-1}\cdot y)|^p \big)\, .
$$

Those coefficients have been introduced in \cite{CDJ}, 
in the setting of products of iid matrices in 
$GL_d(\R)$, 
and proved to be very useful in \cite{CDM} and \cite{CDMP}, 
see also \cite{CDM-deviation}.

\medskip

We shall see that those coefficients decrease exponentially fast to 0, as soon as $\mu$ (almost) admits a moment of order 1, while we obtained only a polynomial speed of convergence in the case of $GL_d(\R)$.

Actually, we will prove the result for the stronger coefficients

$$
\tilde \delta_{p,\infty}(n):=
\E\big(\sup_{x,y\in S^+}|\sigma(Y_n, A_{n-1}\cdot x)-\sigma(Y_n,A_{n-1}\cdot y)|^p \big)\, .
$$

\medskip

\begin{Proposition}\label{delta}
Assume that $\mu$ is strictly contracting and almost 
admits a moment of order $p\ge 1$. Then, there exists $0<a<1$
 such that 
 \begin{equation} \label{expo-coef}
  \delta_{p,\infty}(n)\le \tilde \delta_{p,\infty}(n)
 =O(a^n)\, ,
 \end{equation}  and  
\begin{equation}\label{uniform}
\sup_{x,y\in S^+}\sup_{n\in \N}|\sigma(A_n,x)-\sigma(A_n,y)|\in L^p\, .
\end{equation}
In particular,
\begin{equation}\label{norm-v}
\sup_{n\in \N}|\log \|A_n\|-\log v(A_n)|\in L^p\, .
\end{equation}
\end{Proposition}
\noindent {\bf Proof.} Let $n\in \N$. By item $(iv)$ of Lemma \ref{properties} and item $(vii)$ of Proposition \ref{prop-hennion}, for every $x,\, y\in S^+$, we have 
$$
|\sigma(Y_n, A_{n-1}\cdot x)-\sigma(Y_n,A_{n-1}\cdot y)|
\le (4+2\log L(Y_n))d(A_{n-1}\cdot x, A_{n-1}\cdot y)\le 
(4+2\log L(Y_n))c(A_{n-1})\, .
$$
Let $r\in \N$ be as in 
Definition \ref{contracting}. Then, by item $(vi)$ of 
Proposition \ref{prop-hennion}, there exists 
$\varepsilon>0$ such that
\begin{equation}\label{contraction}
\mu^{*r}(c(g)\le 1-\varepsilon)=:\gamma>0\, .
\end{equation}
Hence, if $m=[(n-1)/r]$, 
$$
\E \big[\big(c(A_{n-1})\big)^p\big]\le \prod_{k=1}^m \E\big[\big(c 
(Y_{kr}\cdots Y_{(k-1)r+1})\big)^p\big]\le \big( \gamma(1-\varepsilon)^p+1-\gamma\big)^{m}\, .
$$
This proves the desired exponential convergence of $(\tilde \delta_{p,\infty}(n))_{n\in \N}$. To conclude the proof, using the cocycle property and the triangle inequality in $L^p$, we infer that 
\begin{align}
\nonumber\E\big[\sup_{x,y\in S^+}\sup_{n\in \N}|\sigma(A_n,x)-\sigma(A_n,y)|^p\big]  &  \le r\E\big[\big(2(2+\log L(Y_1))\, \big)^p\big]\Big(\sum_{
m\ge 0}\big( \gamma(1-\varepsilon)^p+1-\gamma\Big)^{m/p}\Big)^p\\
\label{moment-order-p}  & \qquad  \qquad  = \frac{2^pr\E\big[\big(2+\log L(Y_1)\big)^p\big]}{\Big(1-\big(\gamma(1-\varepsilon)^p+1-\gamma\big)^{1/p}\Big)^{p}}\, .
\end{align}
\hfill $\square$



\medskip


\section{The strong law of large numbers}

\setcounter{equation}{0}

Except the $L^1$-convergences, the results of that section are essentially contained in Hennion's paper \cite{Hennion} (where a more general situation is considered), see his Theorem 2 and its proof.

\medskip

We first recall the version of Kingman's subadditive ergodic theorem relevant to our setting (see \cite[Theorems 1 and 2]{Kingman}). The fact that 
$\lambda_\mu$ in the proposition is 
constant follows from Kolmogorov's $0-1$ law.

\begin{Proposition}[Kingman]\label{kingman}
Assume that $\int_G \big|\log \|g\|\, \big|d\mu (g)<\infty$. Then, 
$(\frac1n \log \|A_n\|)_{n\ge 1}$ converges $\P$-a.s. and in $L^1$ to some constant $\lambda_\mu\in \R$. 
\end{Proposition}

\noindent {\bf Remark.} Using that $\|g\|\ge v(g)$ for every $g\in G^+$, we see that $\log^- \|g\|\le \log ^-v(g)$, where $\log^-(x)=\max
(-\log x,0)$ for every $x>0$. In particular, 
if $\mu$ or $\tilde \mu$ admit a moment of order 1, then, 
$\int_G \big|\log \|g\|\, \big|d\mu (g)<\infty$.

\medskip



We then provide the SLLN for various quantities related to 
$(A_n)_{n\in \N}$ and identify the limit under a stronger assumption.

\begin{Theorem}\label{SLLN}
Assume that  $\mu$ is strictly contracting and that $ \mu$  admits a moment of order 1. Then, for every $x\in S^+$, 
\begin{equation}\label{general-SLLN}
\lim_{n\to +\infty}\frac{\sigma(A_n,x)}n
=\lim_{n\to +\infty} \frac{\log v(A_n)}n
=\lim_{n\to +\infty} \frac{\log \kappa(A_n)}n=\lambda_\mu \qquad \mbox{$\P$-a.s.}\, ,
\end{equation}
where $\lambda_\mu=\int_{G\times S^+}\sigma(g,x)d\mu(g)d\nu(x)$. Moreover, the convergences also hold in $L^1$  and,  we even have
\begin{gather*}\big\|\sup_{x\in S^+}\big|\frac{\sigma(A_n,x)}n-\lambda_\mu\big|\, \big\|_1
\underset{n\to +\infty}\longrightarrow 0
\mbox{ and } \, \sup_{x\in S^+}\big|\frac{\sigma(A_n,x)}n-\lambda_\mu\big|
\underset{n\to +\infty}\longrightarrow 0 \mbox{ $\P$-a.s.}
\end{gather*}
\end{Theorem}

\noindent {\bf Remark.} The $\P$-a.s. and $L^1$ convergence of $(\frac1n\log v(A_n))_{n\in \N}$ when $\int_G|\log v(g)|d\mu(g)<\infty$ (which holds if $\mu$ admits a moment of order 1) follow from Kingman's subadditive ergodic Theorem applied to 
$(-\log v(A_n))_{n\in \N}$. The formula for $\lambda_\mu$ may be 
derived from the formula in the middle of page 1568 of \cite{Hennion}.

\medskip




\noindent {\bf Proof.}  By Proposition \ref{kingman} and the remark after it, we have the $\P$-a.s. and $L^1$ convergence of 
$((\log \|A_n\|)/n)_{n\in \N}$ to $\lambda_\mu$.  

\medskip

By \eqref{norm-v}, we infer the $L^1$ convergence  for $v(A_n)$.

\medskip 

Define $Z:= \sup_{n\in \N}|\log \|A_n\|-\log v(A_n)|$. By \eqref{norm-v}, $Z\in L^1$ and, for every $\varepsilon >0$, 
$$
\sum_{n\in\N}\P(|\log \|A_n\|-\log v(A_n)|\ge \varepsilon n)
\le C \E(Z) <\infty\, .
$$
The $\P$-a.s. convergence for $(v(A_n))_{n\in \N}$ then follows from the one for $(\|A_n\|)_{n\in \N}$ and the Borel-Cantelli lemma. 

\medskip

 The convergences for $\kappa(A_n)$ follows from the bounds $v(A_n)\le \kappa(A_n)\le \|A_n\|$ (see \eqref{v-radius} for the first bound).
 
\medskip 
 
  Finally, notice that 
for every $n\in \N$, 
$$
\sup_{x\in S^+} |\sigma(A_n, x)-n\lambda_\mu|\le \max 
(|\log \|A_n\|- n\lambda_\mu|,|\log v(A_n)-n\lambda_\mu|)\, ,
$$
which proves the remaining convergences.

\medskip

Hence, it remains to identify $\lambda_\mu$. From the above, using the $\mu$-invariance of $\nu$, we infer that 
\begin{align*}
\int_{G\times S^+} \sigma(g,x)d\mu(g)d\nu(x)  &  =\frac1n\int_{S^+}
\E\Big(\sum_{k=1}^n \sigma(Y_k, A_{k-1}\cdot x)\Big) d\nu(x)\\
  &  =\frac1n\int_{S^+}\E(\sigma(A_n, x))d\nu(x)\underset{n\to +\infty}\longrightarrow \lambda_\mu\, .
\end{align*} \hfill $\square$

\medskip

We shall now consider the case of matrix coefficients. The proof 
will relie on Lemma \ref{prop-BDGM} below, which is  essentially Lemma 2.1 of \cite{HH} (see also  Lemma 6.3 of \cite{BDGM} for \eqref{v}). We need also some further notations.

\medskip

For every $0<\delta\le 1$, set 
$$G_\delta:= \{g\in G\, :\, \langle x, gy \rangle \ge \delta\, \, \, \forall x,y\in S^+\}\, ,
$$
and notice that $G^+=\cup_{\delta\in (0,1]}G_\delta$.

\medskip

Let $r\in \N$ be such that $\mu^{*r}(G^+)>0$. There exists $n_0\in \N$, such that $\mu^{*r}(G_{1/n_0})>0$. Then, we define
\begin{equation}\label{Tn}
T_{n_0}:=\inf\{m\in \N\, :\, Y_{mr}\ldots Y_{(m-1)r+1}\in G_{1/n_0}\}\, .
\end{equation}

\medskip

Since $(Y_{mr}\ldots Y_{(m-1)r+1})_{m\in \N}$ is iid with law $\mu^{*r}$ and  $\mu^{*r}(G_{1/n_0})>0$, we know that $T_{n_0}<\infty$ $\P$-a.s.


\begin{Lemma}\label{prop-BDGM}
Assume that $\mu$ is strictly contracting. With the above notations,
\begin{equation}\label{v}
\inf_{n\in \N}\inf_{x\in S^+}\frac{\|A_n x\|}{\|A_n\|}
=\inf_{n\in \N}\frac{v(A_n)}{\|A_n\|}\ge 
\frac{1}{n_0}\min_{1\le n\le rT_{n_0}}\frac{v(A_n)}{\|A_n\|}>0 \qquad 
\mbox{$\P$-a.s.}
\end{equation}
and 
\begin{equation}\label{coeff}
\inf_{n \ge rT_{n_0}}\inf_{x,\, y\in S^+}\frac{\langle y, A_n x\rangle}{
\|Y_1^t\cdots Y_n^t\, y\|}>0 \qquad 
\mbox{$\P$-a.s.}
\end{equation}
\end{Lemma}

\noindent {\bf Proof.}  Let $x\in S^+$. Let $n\in \N$ be such that $n\ge rT_{n_0}$. Using the definition of the action of $G$ on $S^+$ and the definition of $G_{1/n_0}$, we see that 
\begin{align*}
\|A_nx\|  &  = \| Y_n\cdots Y_{rT_{n_0}+1} \big(Y_{rT_{n_0}}\cdots
Y_{r(T_{n_0}-1)+1} \cdot (A_{r(T_{n_0}-1)} x)\big)\|\, \|A_{rT_{n_0}}x\|\\
  &  \ge d\| Y_n\cdots Y_{rT_{n_0}+1}e\|/n_0\, \|A_{rT_{n_0}}x\|\\
  &  \ge \| Y_n\cdots Y_{rT_{n_0}+1}\|\, \|A_{rT_{n_0}}x\|/n_0
  \ge \|A_n\| \frac{\|A_{rT_{n_0}}x\|}{n_0 \|A_{rT_{n_0}}\|}\, ,
\end{align*}
where we used item $(ii)$ of Lemma \ref{properties} for the second inequality.

\medskip

Hence 
$$
\|A_n x\|/\|A_n\|\ge v(A_{rT_{n_0}})/(n_0\|A_{rT_{n_0}}\|) 
{\bf 1}_{\{rT_{n_0}\le n\}}+v(A_n)/\|A_n\|{\bf 1}_{\{rT_{n_0}
> n\}}\, , $$
which proves \eqref{v}.

\medskip

Let us prove \eqref{coeff}. We proceed similarly. Let $x,\, y\in S^+$. Let $n\ge rT_{n_0}$. We have 
\begin{align*}
\langle y, A_n x\rangle  &  =\langle y, Y_n\cdots Y_{rT_{n_0}+1} \big(Y_{rT_{n_0}}\cdots
Y_{r(T_{n_0}-1)+1} \cdot (A_{r(T_{n_0}-1)} x)\big)\rangle \|A_{rT_{n_0}}x\|\\
  &  \ge \|Y_{rT_{n_0}+1}^t\cdots Y_n^ty\|\, \|A_{rT_{n_0}}x\|/n_0\\
  &  \ge \frac1{n_0}\frac{\|Y_{1}^t\cdots Y_n^ty\|\, \|A_{rT_{n_0}}x\|}{\|Y_1^t\cdots Y_{rT_{n_0}}^y\|}=
  \frac1{n_0}\frac{\|Y_{1}^t\cdots Y_n^ty\|\, \|A_{rT_{n_0}}x\|}{\|A_{rT_{n_0}}\|}\, ,
\end{align*}
and \eqref{coeff} follows from \eqref{v}. \hfill $\square$

\bigskip

We denote by $\tilde \mu$ the pushforward image of $\mu$ by the map $g\to g^t$.

\begin{Theorem}\label{SLLN-coeff}
Assume that  $\mu$ is strictly contracting and that
$\tilde \mu$ admits a moment of order 1. Then, 
\begin{equation*}\left(\sup_{x,\, y\in S^+}\left|
\frac{\log \langle y, A_n x\rangle}n-\lambda_\mu\right|
\right)_{n\in \N}\underset{n\to +\infty}\longrightarrow 0\, \mbox{ $\P$-a.s.} 
\end{equation*} In particular, 
\begin{equation*}\left(\left|
\frac{\inf_{x,\, y\in S^+}\log \langle y, A_n x\rangle}n-\lambda_\mu\right|
\right)_{n\in \N}\underset{n\to +\infty}\longrightarrow 0 
\, \mbox{ $\P$-a.s.}
\end{equation*} 

\noindent Moreover, 
$((\log \|A_n\|-n\lambda_\mu)/n)_{n\in \N}$ and 
$((\log \kappa(A_n)-n\lambda_\mu)/n)_{n\in \N}$ converge $\P$-a.s. and in $L^1$ to 0; and $((\log v(A_n)-n\lambda_\mu)/n)_{n\in \N}$  converges $\P$-a.s. to 0.
\end{Theorem}


\noindent {\bf Proof.} First notice that Proposition \ref{kingman} applies, which yields the $\P$-a.s. and 
$L^1$ convergence for $\log \|A_n\|$ and for $\log \|A_n^t\|$ by item $(iii)$ of Lemma \ref{properties}. 

\medskip

By Lemma  \ref{prop-BDGM}, there exists a random variable $W\ge 0$ such that, for every $x,\, y\in S^+$ and every $n\in \N$, 
 on the set  $\{rT_{n_0}\le n\}$,
\begin{equation}\label{coeff-SLLN}
0\le 
\log \|A_n\|- \log \langle y, A_n x\rangle \le \log W+
\log\|A_n\|- \log \|Y_1^t \cdots Y_n^ty\|\, .
\end{equation}

\medskip

Let $\varepsilon>0$. Using that $(Y_1, \ldots , Y_n)$ and $(Y_n, \ldots , Y_1)$ have the same law, we get  

\begin{gather*}
\sum_{n\ge 1}\P(\sup_{y\in S^+} \big| \log \|Y_1^t \cdots Y_n^ty\|-
\log \|Y_1^t \cdots Y_n^te\|\, \big|\ge \varepsilon n)\\
\le \sum_{n\ge 1} \P(\sup_{y\in S^+} \sup_{m\in \N}\big| \log \|Y_m^t \cdots Y_1^ty\|-
\log \|Y_m^t \cdots Y_1^te\|\, \big|\ge \varepsilon n)
<\infty\, ,
\end{gather*}
where we used Proposition \ref{delta} for $\tilde \mu$.

By the Borel-Cantelli lemma, using item $(ii)$ of Lemma \ref{properties}, we infer that 
$$
\frac{\sup_{y\in S^+} \big| \log \|Y_1^t \cdots Y_n^ty\|-
\log \|A_n^t\|\, \big|}n\underset{n\to +\infty}
\longrightarrow 0 \mbox{\quad $\P$-a.s.}
$$
Combining this with \eqref{coeff-SLLN} (recall that 
$\P(T_{n_0}<\infty)=1$ and that $\|g\|
\le d\|g^t\|$ for every $g\in G$) we obtain that 
$$
\sup_{x,\, y\in S^+}\frac{\big|\log \|A_n\|- \log \langle y, A_n x\rangle\big|}{n}\underset{n\to +\infty} 
\longrightarrow 0\quad \mbox{$\P$-a.s.}\, 
$$
This gives the desired convergence for the coefficients. 
The $\P$-a.s. convergences for $\kappa(A_n)$ and $v(A_n)$ 
follow from the inequalities
$$
\frac{\inf_{x,\, y\in S^+}\log \langle y, A_n x\rangle}n
\le \frac{\log v(A_n)}n\le \frac{\log \kappa(A_n)}n 
\le \frac{\log \|A_n\|}n\, .
$$
The $L^1$ convergence for $\kappa(A_n)$, follows from 
Theorem \ref{SLLN} applied to $\tilde \mu$, using 
item $(iii)$ of Lemma \ref{properties}, noticing that 
$(Y_1, \ldots , Y_n)$ has the same law as $(Y_n, \ldots , Y_1)$.  \hfill $\square$

\medskip

Assume that  $\mu$ (hence $\tilde \mu$) is strictly contracting and that  $\mu$  and 
$\tilde \mu$ both admit a moment of order 1. Denoting by $\tilde \nu$   
the only $\tilde \mu$-invariant probability on $S^+$, and using 
that $A_n^t$ and $Y_n^t\ldots Y_1^t$ have the same law, we have 
\begin{align*}
\int_{G\times S^+}\sigma(g,x)d\mu (g)d\nu (x)  &  =\lambda_\mu 
=\lim_{n\to +\infty}\frac{\E[\log \|A_n\|]}n \\
&  = \lim_{n\to +\infty}\frac{\E[\log \|A_n^t\|]}n
  =\lim_{n\to +\infty}\frac{\E[\log \|Y_n^t\ldots Y_1^t\|]}n\\
&  =\lambda_{\tilde \mu}=\int_{G\times S^+}\sigma(g,x)d\tilde \mu (g)d\tilde\nu (x)
\, .
\end{align*}

\bigskip

Under our assumptions, one cannot expect the $L^1$ convergence in Theorem \ref{SLLN-coeff} for $v(A_n)$. 

\medskip

For instance take $\mu $ such that for every $n\in \N$, 
$\mu(\{g_n\})= \frac1{\pi^2n^2}$ and $\mu (\{h\})=5/6$, 
with $g_n=\left(\begin{array}{cc}
1  &  0\\
0  &  2^{-n}
\end{array}\right)$ and $h=\left(\begin{array}{cc}
1  &  1\\
1  &  1
\end{array}\right)$. Then, for any $k_1, \ldots , k_r \in \N$, $v(g_{k_1}\cdots g_{k_r})\le v(g_{k_r})\le 2^{-k_r}$.

\medskip Hence $\E(\log v(A_n))\le \frac1{6^{n-1}} \sum_{
k\in \N} \frac{-k}{\pi^2k^2}=-\infty$.

\medskip

Similarly, even if $\mu$ and $\tilde \mu$ are strictly contracting and admit a moment of order 1, we may not have 
$L^1$ convergence for the coefficients.

\medskip

For instance, let 
$\mu$ be such that $\mu(\{Id\})=1/2$, with $Id$ the 
identity  matrix. Then, $\mu^{*n}(\{Id\})\ge 2^{-n}$ and, with $\{e_1,e_2\}$  the canonical basis of $\R^2$, 
$\mu (\{g\in G\, :\, \langle e_1,g e_2\rangle=0\})>0$, so that $\E(\log \langle e_1, A_n e_2\rangle) =-\infty$.

\section{The CLT and the asymptotic variance}

\setcounter{equation}{0}

We start by proving a martingale-coboundary decomposition. In the case of invertible matrices, such a decomposition was only available for $p\ge 2$ while here it holds as soon as $p\ge 1$.

\begin{Proposition}\label{martingale-decomposition}
Assume that  $\mu$ is strictly contracting and  admits a moment of order $p\ge 1$. There exists a continuous and bounded function 
$\psi$ on $X$ such that
$\big( \sigma(Y_n, A_{n-1}\cdot x)-\lambda_\mu+\psi(A_n\cdot x)-\psi(A_{n-1}
\cdot x)\big)_{n\in \N}$ is a sequence of martingale differences in $L^p$. 
If moreover $W_0$ is a random variable with law $\nu$, independent from 
$(Y_n)_{n\in \N}$, then $\big( \sigma(Y_n, A_{n-1}\cdot W_0)-\lambda_\mu+\psi(A_n\cdot W_0)-\psi(A_{n-1}
\cdot W_0)\big)_{n\in \N}$ is a \emph{stationary} and ergodic sequence of martingale differences in $L^p$. 
\end{Proposition}

\noindent {\bf Remark.} The function $\psi$ in the theorem 
is given by 
\begin{equation}\label{psi}
\psi(x):= 
\sum_{n\ge 1}\Big(\int_{G\times G}\sigma(g, g'\cdot x)d\mu(g)d\mu^{*(n-1)}(g')-\lambda_\mu\Big)\, .
\end{equation}  

\medskip

\noindent {\bf Proof.} Let   $\psi$ be given by \eqref{psi}. The fact that $\psi$ is  well-defined and continuous follows from Proposition \ref{delta}.

\medskip

 Then, notice that 
$$
\sigma(g,x)-\lambda_\mu= \sigma(g,x)-\int_G\sigma(g',x)d\mu(g')+ 
\int_G\sigma(g',x)d\mu(g')-\lambda_\mu 
$$
and, using the definition of $\psi$,
$$
\int_G\sigma(g,x)d\mu(g)-\lambda_\mu+\int_G\psi(g\cdot x)d\mu(g)
=\psi(x)\, .
$$
Now, $\Big(\sigma(Y_n,A_{n-1}\cdot x)-\int_G\sigma(g,A_{n-1}\cdot x)d\mu(g)\Big)_{n\in \N}$ is a sequence of martingale differences in $L^p$ (notice that 
$x\mapsto \int_G\sigma(g,x)d\mu(g)$ is bounded). Moreover, 
$$
\int_G\sigma(g,A_{n-1}\cdot x)d\mu(g)-\lambda_\mu+\psi(A_n \cdot x)- 
\psi(A_{n-1}\cdot x)=\psi(A_n\cdot x)-\int_G\psi(g A_{n-1}\cdot x)\,d\mu(g) ,
$$
and the RHS defines a sequence of bounded martingale differences. 

\medskip

 The final statement follows from the fact that 
$((Y_n,A_{n-1}\cdot W_0))_{n\in \N}$ is a stationary and (uniquely) ergodic Markov chain.

\hfill  $\square$

\medskip



\begin{Definition}\label{aperiodic}
We say that a probability $\mu$ on $G$ is aperiodic if 
the group generated by $\{\log \kappa(g)\, :\, g\in \Gamma_\mu\}$ is dense in $\R$.
\end{Definition}

We now state and prove various CLTs. Those CLTs are proved in Hennion \cite{Hennion} by a slightly different approach (also based on a martingale-coboundary 
decomposition). We complement the results of Hennion by identifying the asymptotic variance $s^2$ in several ways and by characterizing the fact that $s^2>0$. The characterization is the same as in \cite{BDGM} or in \cite{GLX} but its proof does not require exponential moments as in those works.  

\begin{Proposition}\label{prop-CLT}
Assume that $\mu$ is strictly contracting  and admits a moment of order $2$. Then, there exists $s^2\ge 0$ such that, with $W_0$ as in Proposition \ref{martingale-decomposition}, $\frac1n\E[(\sigma(A_n,W_0)-n \lambda_\mu)^2]\underset{n\to +\infty}\longrightarrow s^2$ and $\frac1{\sqrt n} (\sigma(A_n,W_0)-n \lambda_\mu)
\Rightarrow {\mathcal N}(0,s^2)$. If there does not exist $m\in \N$ and 
$\psi_m$ continuous on $S^+$ such that 
\begin{equation}\label{null-variance}
\sigma(g,x)-m\lambda_\mu =\psi_m(x)-\psi_m(g\cdot x)\qquad \mbox{ for $\mu^{\otimes m}\otimes \nu$-almost every $(g,x)\in G\times S^+$}\, ,
\end{equation}
then $s^2>0$. In particular, if $\mu$ is aperiodic, then 
$s^2>0$.
\end{Proposition}
\noindent {\bf Remark.} Under the assumptions of the 
proposition we actually have the functional central limit theorem. It is well-known  that the  variance is given by  
\begin{align*}
s^2  &  =\E(\sigma(A_1,W_0)^2)+2\sum_{n\ge 2}\E(\sigma(A_1,W_0)\sigma(A_n, W_0))\\
  &  =\int_{G\times S^+}\sigma^2(g,x)d\mu(g)d\nu(x)+2\sum_{n\ge 2} 
\int_{G^2\times S^+}\sigma(g,x)\sigma(g'g,x)d\mu^{*(n-1)}(g')d\mu(g)d\nu(x)\, .
\end{align*}

\medskip

\noindent {\bf Proof.} For every $n\in \N$, set 
$D_n:=\sigma(Y_n, A_{n-1}\cdot W_0)-\lambda_\mu+\psi(A_n\cdot W_0)-\psi(A_{n-1}
\cdot W_0)$. By Proposition \ref{martingale-decomposition}, $(D_n)_{n\in \N}$ is a stationary and ergodic sequence of martingale differences in $L^2$. In particular, 
$(D_1+\ldots +D_n)/\sqrt n\Rightarrow {\mathcal N}(0,s^2)$, with 
$s^2=\E(D_1^2)=\E((D_1+\ldots +D_n)^2)/n$.  Hence, the CLT with the description of the variance follows from the following reformulation of Proposition \ref{martingale-decomposition}:
\begin{equation}\label{martingale-coboundary}
\sigma(A_n, W_0)-n\lambda_\mu = (D_1+\ldots +D_n)+\psi(W_0)-\psi (A_n\cdot W_0)\, .
\end{equation}
Assume now that $s^2=0$. Then 
$$
\int_G (\sigma(g,x)-\lambda_\mu -\psi(x)+\psi(g\cdot x))^2\,  d\mu(g)d\nu(x)=0\, . 
$$
Hence, \eqref{null-variance} holds with $m=1$ and $\psi_1=\psi$. Let $m>1$. Notice that $\mu^{*m}$  is strictly contracting and
 admits a moment of order $p$  and that the unique 
$\mu^{*m}$-invariant measure is the unique $\mu$-invariant measure. Notice also that $\lambda_{\mu^{*m}}=m\lambda_\mu$. Applying the above argument to $\mu^{*m}$, we infer that 
there exists a continuous $\psi_m$ satisfying to \eqref{null-variance}. 

\medskip

Using that $\psi_m$ is continuous, we see that \eqref{null-variance} holds for every $g\in {\rm supp}\, \mu^{*m}$ and every $x\in {\rm supp}\, \nu$.

\medskip
Let $g\in  {\rm supp}\, \mu^{*m}\subset \Gamma_\mu$. Then, 
$v_g\in \Lambda_\mu\subset {\rm supp}\, \nu$. Since
$g\cdot v_g=v_g$ and $\sigma(g,v_g)=\log \kappa(g)$, we infer that $\psi_m(g\cdot v_g)=\psi_m(v_g)$ and that  
$\log \kappa(g)= m\lambda_\mu$. 

\medskip

Hence, $\log \kappa(\Gamma_\mu)\subset \lambda_\mu \N$ and 
$\mu$ cannot be aperiodic. \hfill $\square$

\medskip

\begin{Proposition}\label{variance-identification}
Assume that $\mu$ is strictly contracting  and  admits a moment of order  2. Then, 
with $s^2$ as in Proposition \ref{prop-CLT},
\begin{align*}
s^2  &  =\lim_{n\to +\infty}\frac1n \sup_{x\in S^+}\E(
(\sigma(A_n, x)-n\lambda_\mu)^2)\\
  &  =\lim_{n\to +\infty}\frac1n \E(
(\log \|A_n\|-n\lambda_\mu)^2)\\
  &  =\lim_{n\to +\infty}\frac1n \E(
(\log v(A_n)-n\lambda_\mu)^2)\\
  & = \lim_{n\to +\infty}\frac1n \E(
(\log \kappa(A_n)-n\lambda_\mu)^2)\, ,
\end{align*}
and the CLT holds if we replace $\sigma(A_n,W_0)$ with 
$\sigma(A_n,x)$, $\log \|A_n\|$, $\log v(A_n)$ or $\log \kappa(A_n)$. Moreover 
$$
\sup_{x\in S^+}\sup_{t\in \R}\Big|\P( \sigma(A_n, x)-n
\lambda_\mu\le t\sqrt n)-\phi(t/s^2)\Big|\underset{n\to +\infty}\longrightarrow 0\, .
$$
\end{Proposition}
\noindent {\bf Proof.} The result follows  from Proposition 
\ref{prop-CLT} and  Proposition 
\ref{delta} (using \eqref{v-radius}). \hfill $\square$

\medskip

We also have a (functional)  CLT for the coefficients. As noticed in the previous section, one cannot expect in general 
to identify $s^2$ thanks to the matrix coefficients as in Proposition \ref{variance-identification}.

\medskip

\begin{Proposition}
Assume that $\tilde \mu$ is strictly contracting and  admits a moment of order 2. Then, 
\begin{gather*}
\sup_{x,\, y\in S^+}\sup_{t\in \R}\Big|\P(\log \langle x, A_n y
\rangle )-n
\lambda_\mu\le t\sqrt n)-\phi(t/s^2)\Big|\underset{n\to +\infty}\longrightarrow 0\, ,\\
\sup_{t\in \R}\Big|\P(\inf_{x,\, y\in S^+}\log \langle x, A_n y
\rangle )-n
\lambda_\mu\le t\sqrt n)-\phi(t/s^2)\Big|\underset{n\to +\infty}\longrightarrow 0\,
\end{gather*}
\end{Proposition}

\noindent {\bf Proof.} We proceed as for the proof 
of Theorem \ref{SLLN-coeff}. By Proposition \ref{delta} 
applied with $\tilde \mu$, 
\begin{gather*}
\sum_{n\in \N}\P(\sup_{y\in S^+}\big|\log \|Y_1^t\cdots Y_n^t\|- \log \|Y_1^t\cdots Y_n^ty\|\, \big|\ge \varepsilon \sqrt n)
\\\le \sum_{n\in \N} \P(\sup_{y\in S^+}\sup_{m\in \N}\big|\log \|Y_m^t\cdots Y_1^t\|- \log \|Y_m^t\cdots Y_1^ty\|\, \big|\ge \varepsilon \sqrt n)<\infty\, .
\end{gather*}
Using \eqref{coeff-SLLN}, the fact that $\P(T_{n_0}<\infty)=1$,  and Proposition 
\ref{variance-identification} with $\tilde \mu$, the result follows. \hfill $\square$

\section{The almost sure invariance principle}

\setcounter{equation}{0}

\begin{Theorem}\label{ASIP}
Let $p\ge 2$. Assume that $\mu$ is strictly contracting  and  admits a moment of order $p$. Let $s^2$ be as in Proposition 
\ref{prop-CLT}. Then, one can redefine the process 
$(\sigma(A_n, W_0))_{n\in \N}$ on another probability space on which there exist iid  variables $(N_n)_{n\in \N}$ with law ${\mathcal N}(0,s^2)$, such that
\begin{gather*}
|\sigma(A_n,W_0)- n\lambda_\mu-(N_1+\ldots +N_n)|= o(\sqrt{n\log\log n})\quad \mbox{$\P$-a.s. if $p=2$}\\
\mbox{and}\quad |\sigma(A_n,W_0)- n\lambda_\mu-(N_1+\ldots +N_n)|= o(n^{1/p})\quad \mbox{$\P$-a.s. if $p>2$}
\end{gather*}
\end{Theorem}

\noindent {\bf Remark.} It is not necessary here that
 $s^2>0$.

 \noindent {\bf Proof.} When $p>2$, the result follows from  Theorem 1 of \cite{CDM} by taking  into account \eqref{expo-coef}. The case $p=2$ follows from \eqref{martingale-coboundary} and the ASIP for martingales with stationary and ergodic increments in $L^2$.  \hfill $\square$
 
 \medskip

 Proceeding as in the proof of Theorem \ref{SLLN}, one can prove that the above theorem holds if we replace 
 $(\sigma(A_n,W_0))_{n\in \N}$ with any of the following sequences: $(\sigma(A_n,x))_{n\in \N}$ (for a given $x\in S^+$), $(\log\|A_n\|)_{n\in \N}$, $(\log \kappa(A_n))_{n\in \N}$ or $(\log v(A_n))_{n\in \N}$.

 \medskip
 

Let us give the ASIP for the matrix coefficients.

\medskip

 \begin{Theorem}
 Let $p\ge 2$. Assume that $\mu$ is strictly contracting  and  that $\mu$ and $\tilde \mu$ admit a moment of order $p$. Then, for every $x,\, y\in S^+$, one can redefine the process 
$(\log \langle y, A_n x\rangle)_{n\in \N}$ on another probability space on which there exist iid  variables $(N_n)_{n\in \N}$ with law
${\mathcal N}(0,s^2)$, such that
\begin{gather*}
|\log \langle y, A_n x\rangle- n\lambda_\mu-(N_1+\ldots +N_n)|= o(\sqrt{n\log\log n})\quad \mbox{$\P$-a.s. if $p=2$}\\
\mbox{and}\quad|\log \langle y, A_n x\rangle- n\lambda_\mu-(N_1+\ldots +N_n)|= o(n^{1/p})\quad \mbox{$\P$-a.s. if $p>2$}
\end{gather*}
 \end{Theorem}
 
 The proof may be done similarly to the one of Theorem \ref{SLLN-coeff}. Since $\tilde \mu$ almost admit a moment of order $p\ge 1$, 
 $$
 \frac{\sup_{y\in S^+}\Big|\log \|Y_1^t\cdots Y_n^t\|- 
 \log \|Y_1^t\cdots Y_n^ty\|\, \Big|}{n^{1/p}}\underset{n\to +\infty}\longrightarrow 0\qquad \mbox{$\P$-a.s.}\, ,
 $$
 and we conclude thanks to Theorem \ref{ASIP}, using 
\eqref{coeff-SLLN} and the fact that $\P(T_{n_0}<\infty)=1$.

\medskip

In case of exponential moments, combining ideas from \cite{CDM} and \cite{CDKM}, it is possible to obtain logarithmic rates in the ASIP. However those rates are not as good as for the sums of independent variables: in the case of a sum of iid variables it is possible to obtain a rate $O((\log n)^(1/_gamma))$ instead of $O((\log n)^{2+
1/ \gamma})$ under exponential moments of order $\gamma\in (0,1]$. Let us state the results, the proof will be done in a forthcoming work \cite{CDM-asip}. 

\medskip

\begin{Theorem}\label{ASIP-expo}
 Assume that $\mu$ is strictly contracting  and  admits an exponential moment of order $\gamma\in (0,1]$. Let $s^2$ be as in Proposition 
\ref{prop-CLT}. Then, one can redefine the process 
$(\sigma(A_n, W_0))_{n\in \N}$ on another probability space on which there exist iid variables $(N_n)_{n\in \N}$ with law ${\mathcal N}(0,s^2)$,  such that
\begin{gather*}
|\sigma(A_n,W_0)- n\lambda_\mu-(N_1+\ldots +N_n)|= 
O((\log n)^{2+1/\gamma})\quad \mbox{$\P$-a.s.}
\end{gather*}
\end{Theorem}

Again, the theorem is true if if we replace 
 $(\sigma(A_n,W_0))_{n\in \N}$ with any of the following sequences: $(\sigma(A_n,x))_{n\in \N}$, $(\log\|A_n\|)_{n\in \N}$, $(\log \kappa(A_n))_{n\in \N}$ or $(\log v(A_n))_{n\in \N}$.
 
 \medskip
 
 We also have a result for the coefficients.

\begin{Theorem}\label{ASIP-expo-coeff}
 Assume that $\mu$ is strictly contracting  and that $\mu$ and $\tilde \mu$  admit an exponential moment of order $\gamma\in (0,1]$. Let $s^2$ be as in Proposition 
\ref{prop-CLT}. Then, for every $x,\, y\in S^+$, one can redefine the process 
$(\log \langle y, A_n x\rangle)_{n\in \N}$ on another probability space on which there exists iid normal variables $(N_n)_{n\in \N}$ 
with law ${\mathcal N}(0,s^2)$ such that
\begin{gather*}
|\log \langle y, A_nx\rangle- n\lambda_\mu-(N_1+\ldots +N_n)|= 
O((\log n)^{2+1/\gamma})\quad \mbox{$\P$-a.s.}
\end{gather*}
\end{Theorem}

\noindent {\bf Proof.} Let $x,y\in S^+$. Let $n\in \N$. We have 
$$
\log \langle y, A_nx\rangle - \log \|A_nx\|= 
\log \langle y, A_n\cdot x\rangle\, .
$$
In view of Theorem \ref{ASIP-expo}, it suffices to prove that 
there exists $c>0$ such that 
$$\sum_{n\ge 1} \P(|\log \langle y, A_n\cdot x\rangle|\ge  
c(\log n)^{1/\gamma})<\infty\,.$$

\medskip

We will use the following simple observation, which follows from the independence of $(Y_n)_{n\in \N}$. For every $x,y\in S^+$, every integers $1\le m\le n$ and every $t>0$ 
$$
\P(|\log \langle y, A_n\cdot x\rangle|\ge t)\le \sup_{u,v\in S^+} \P(|\log \langle u, A_m\cdot v \rangle|\ge t)\, .
$$

\medskip

Let $\eta, \delta$ be as in \eqref{regularity-expo}. 
For $n\ge [(\log n)^{c\gamma/\eta}]$ (with $[\cdot]$ the integer part), using \eqref{regularity-expo}, we have 
\begin{gather*}\P(|\log \langle y, A_n\cdot x\rangle|\ge  
c(\log n)^{1/\gamma})
=\P(|\log \langle y, A_n\cdot x\rangle|\ge  
\eta (\log (n^{(c/\eta)^\gamma}))^{1/\gamma})\\
\le \P\big(\sup_{u,v\in S^+} |\log \langle u, A_{[(\log (n^{(c/\eta)^\gamma}))^{1/\gamma}]}\cdot v\rangle|\ge  \eta[(\log (n^{c\gamma/\eta}))^{1/\gamma}\big)\\
=  o({\rm exp}( -\delta [(\log (n^{(c/\eta)^\gamma}))^{1/\gamma}]^\gamma)=o(n^{-\delta(c/\eta)^\gamma})\, ,
\end{gather*}
and the result follows by taking $c=\eta(2/\delta)^{1/\gamma}$. \hfill $\square$
 
 \section{The Berry-Esseen theorem}
 
 \setcounter{equation}{0}
 
 \subsection{Berry-Esseen for the norm cocycle and the matrix norm}

 \begin{Theorem} \label{berry-esseen}
Let $p\in (2,3]$. Assume that $\mu$ is strictly contracting  and  admits a moment of order $p$. Assume that $s^2>0$ with $s^2$ as in Proposition \ref{prop-CLT}. Then, setting $ \displaystyle v_n =  \Big (  \frac{1}{n }  \Big )^{p/2-1}$, we have 
\beq \label{ineBE1}
\sup_{t \in {\mathbb R}} \Big | {\mathbb P} \big ( \sigma(A_n, W_0)-n\lambda_\mu  \leq t \sqrt{n} \big ) - \Phi (t/ s)  \Big |= O(v_n) \, , \eeq
\beq \label{BE-cocycle}
 \sup_{x \in S^+}\sup_{t \in {\mathbb R}} \Big | {\mathbb P} \big ( \sigma(A_n, x)  - n \lambda_{\mu}  \leq t \sqrt{n} \big ) - \Phi (t/ s)  \Big |= O( v_n) \, ,\eeq
 \beq \label{ineBE2}
\sup_{t \in {\mathbb R}} \Big | {\mathbb P} \big ( \log \|A_n\|-n\lambda_\mu  \leq t \sqrt{n} \big ) - \Phi (t/ s)  \Big | =O( v_n) \, , \eeq
\end{Theorem}

\medskip

\noindent {\bf Proof.}  Redo the proof of Theorem 2.1 of \cite{CDMP} with $T=n^{p/2-1}$, using \eqref{expo-coef}. \hfill $\square$

\medskip

\noindent {\bf Remarks.} By some arguments already mentionned, \eqref{ineBE2} also holds if $\tilde \mu$ 
is  strictly contracting and admits a moment of order $p\in (2,3]$.
Let us notice that  \eqref{ineBE1} follows also from 
Theorem 2.3 of \cite{Jirak-2}, since the Assumptions 2.1 there are satisfied due to the exponential convergence of the coefficients $\delta_{\infty,p}$ in Proposition \ref{delta}.  

Finally, let us mention that  Grama et al. \cite{GLX} obtained  \eqref{BE-cocycle} and \eqref{ineBE2} for $p=3$ under their condition $A.2$. That condition is equivalent to the condition used in 
Theorem \ref{theo-strong-cond} below.

\medskip


\subsection{Berry-Esseen for the spectral radius and the matrix coefficients}






\medskip

\begin{Proposition}\label{BE-radius-polynomial}
Let $p\in (2,3)$. Assume that $\mu$ is strictly contracting,  admits a moment of order $p$ and almost admits a moment 
of order $q\in [p, \max(p,(p-2)/(3-p))]$. Assume that $s^2>0$. Set $ \displaystyle v_n =  \Big (  \frac{1}{n }  
\Big )^{p/2-1}$ if $p\in (2, 1+\sqrt 3]$ and $ \displaystyle v_n =  \Big (  \frac{1}{n }  
\Big )^{q/2(q+1)}$ if $p\in (1+\sqrt 3,3]$. Then,  
\begin{equation}\label{BE-radius-bis}
\sup_{t \in {\mathbb R}} \Big | {\mathbb P} \big ( \log v(A_n)-n\lambda_\mu  \leq t \sqrt{n} \big ) - \Phi (t/ s)  \Big | =O( v_n) \,  
\end{equation}
and
\beq \label{BE-radius}
\sup_{t \in {\mathbb R}} \Big | {\mathbb P} \big ( \log \kappa(A_n)-n\lambda_\mu  \leq t \sqrt{n} \big ) - \Phi (t/ s)  \Big | =O( v_n) \, . \eeq
\end{Proposition}

\noindent {\bf Remark.} When $p\le 1+\sqrt 3$ the condition on $q$ reads $q=p$ hence is satisfied. When $p=3$ the condition on $q$ reads $q\ge p$. \eqref{BE-radius-bis} and \eqref{BE-radius} also hold if $\tilde \mu$ satisfies the assumptions of the proposition. 

\noindent {\bf Proof.} Applying Proposition \ref{delta} (with $p=q$) and Theorem 
\ref{berry-esseen}, we see that we can use Lemma 2.1 of \cite{CDMP} with $T_n=\log \|A_n\|-n \lambda_\mu$, 
$R_n=\log v(A_n)-\log\|A_n\|$, $a_n=n^{(p-2)/2}$, 
$b_n=n^{q/2(q+1)}$ and $c_n=(\sqrt n/b_n)^q)$ to obtain 
\eqref{BE-radius-bis}. Then, \eqref{BE-radius} follows from the fact that $v(A_n)\le \kappa(A_n)\le \|A_n\|$. \hfill $\square$

\begin{Proposition}\label{BE-radius-expo}
Assume that $\mu$ is strictly contracting,  admits a moment of order $3$ and almost admits an exponential moment of order 
$\gamma\in (0,1]$. Assume that $s^2>0$. Set $ \displaystyle v_n =  \frac{(\log n)^{1/\gamma}}{n ^{1/2}}$. Then,  
\begin{equation*}
\sup_{t \in {\mathbb R}} \Big | {\mathbb P} \big ( \log v(A_n)-n\lambda_\mu  \leq t \sqrt{n} \big ) - \Phi (t/ s)  \Big | =O( v_n) \,  
\end{equation*}
and
\beq \label{BE-radius-ter}
\sup_{t \in {\mathbb R}} \Big | {\mathbb P} \big ( \log \kappa(A_n)-n\lambda_\mu  \leq t \sqrt{n} \big ) - \Phi (t/ s)  \Big | =0( v_n) \, . \eeq
\end{Proposition}
\noindent {\bf Remarks.} \eqref{BE-radius-ter} also holds if $\tilde \mu$ satifies the assumptions of the proposition. \eqref{BE-radius-ter} has been proved in \cite{GLX} under a much stronger assumption than exponential moments.

\noindent {\bf Proof.} As before we prove the result for 
$v(A_n)$ in place of $\kappa(A_n)$.  
\medskip
Let $\varepsilon\in (0,1)$ be such \eqref{contraction} holds.
  Let $x, \,y\in S^+$. Let $n\in \N$. Let $\omega\in \Omega$. 
Let $1\le m< [n/r]$ be such that 
$c(Y_{mr}\cdots Y_{(m-1)r+1})(\omega)\le 1-\varepsilon$. Using the cocycle property and several items  of Proposition \ref{prop-hennion} 
(in particular item $(iv)$), we see that 
\begin{align*}
|\sigma(A_n, x)-\sigma(A_n,y)|  &  \le |
\sigma(Y_n\cdots Y_{mr+1}, A_{mr}\cdot x)-\sigma(
Y_n\cdots Y_{mr+1}, A_{mr}\cdot y)| +|\sigma(A_{mr}, x)-\sigma(A_{mr},y)|\\ 
&  \le 2 \ln \big(1/(1-d(A_{mr}\cdot x,A_{mr}\cdot y))\big)+\log \|A_{mr}\|- \log v(A_{mr})\\
  &  \le 2\ln(1/\varepsilon)+\log \|A_{mr}\|- 
  \log v(A_{mr})\, .
\end{align*}
Define
\begin{equation}\label{gamma_m}\Gamma_m:=\{\exists k\in {1,\ldots, m}\, :\, 
c(Y_{kr}\cdots Y_{(k-1)r+1})\le 1-\varepsilon\}\,.
\end{equation} 
Taking the supremum over $x$ and the infimum over $y$, we infer that on $\Gamma_m$,
$$
\log \|A_{n}\|- 
  \log v(A_{n})\le 2\ln(1/\varepsilon)+
\max_{1\le k\le m} \big(\log \|A_{kr}\|-  \log v(A_{kr})\big)\, .
$$
Hence, for $\eta m\ge \ln(1/\varepsilon)$, 
using Lemma \ref{deviation} below,
\begin{gather*}
\P(\log \|A_{n}\|- \log v(A_{n})\ge 2\eta m)\le 
\P(\Gamma_m^c)+ \P\big(\max_{1\le k\le m} (\log \|A_{kr}\|-  \log v(A_{kr}))\ge \eta m\big)\\
\le (1-\gamma)^m +C_\eta{\rm e}^{-\delta_\eta n}\, .
\end{gather*}
Taking $m\sim C\log n$, with $C|\log (1-\gamma)|> 1/2$, we infer that the right-hand side is bounded by $D/\sqrt n$, and we conclude thanks to Lemma 2.1 of \cite{CDMP}, using Theorem \ref{berry-esseen}. \hfill $\square$

\medskip

\begin{Lemma}\label{deviation}
Assume that $\mu$ is strictly contracting and almost admits some exponential moment of order $\gamma\in (0,1]$. Then, there exist $ \eta, \, \delta>0$ 
 such that 
$$
\P(\max_{1\le k\le n} \big|\log v(A_k)- \log \|A_k\|\, \big|\ge \eta n)\le 
{\rm e}^{-\delta n^\gamma }\, .
$$
\end{Lemma}

\noindent {\bf Proof.} For every $n\in \N$, using that 
$\|\cdot \|$ is submultiplicative and that $v$ is supermultiplicative, we see that, setting $\tau:=\E(\log\|Y_1\|/v(Y_1))$, 
$$
\max_{1\le k\le n} \big|(\log (\|A_k\|)- \log (v(A_k))\, \big|
\le \max_{1\le k\le n}\Big|\sum_{i=1}^k 
\big[\log\big(\|Y_i\|/v(Y_i)\big)-\tau\big]\Big|+n\tau\, .
$$ 
Then the desired result follows from Theorem 2.1 of 
\cite{FGL-deviation}, see their estimate $(2.7)$.\hfill $\square$

\medskip

\begin{Proposition}
Let $p\in (2,3]$. Assume that $\mu$ or $\tilde \mu$ satisfies the assumptions of Proposition 
\ref{BE-radius-polynomial} if $p<3$ and those of Proposition 
\ref{BE-radius-expo} if $p=3$. Then, 
\eqref{BE-radius} (if $p<3$) and \eqref{BE-radius-ter} 
(if $p=3$) hold with $\inf_{x,y\in S^+}\langle y , A_n x\rangle$ in place of $\kappa(A_n)$.
\end{Proposition}

\noindent {\bf Proof.} Recall that, for every $0< \delta \le 1$, we defined 
$$
G_{\delta}:=\{g\in G\, :\, \langle x,g\cdot y\rangle \ge \delta \quad \forall x,y\in S^+\}\, .
$$

Notice that $g\in G_\delta$ if and only if for every 
$y\in S^+$ all the coordinates of $g\cdot y$ are greater that $\delta$, i.e. $g\cdot y -d\delta e \in (\R^+)^d$. 

\medskip
Let $x, \,y\in S^+$. Let $n\in \N$, $\omega\in \Omega$ and $
n_0\in \N$. 
Let $1\le m< [n/r]$ be such that 
$(Y_{mr}\cdots Y_{(m-1)r+1})(\omega)\in G_{1/n_0}$. We have 
(omitting $\omega$) 
\begin{gather*}
\langle y, A_n x\rangle \ge \langle Y_{mr+1}^t
\cdots Y_n^ty, \frac{A_{mr}x}{\|A_{mr}x\|}\rangle \|A_{mr}x\|\ge (1/n_0) \|Y_{mr+1}^t \cdots Y_n^ty\|\, \frac{\|A_nx\|}
{\|Y_n \cdots Y_{mr+1}\|}\, .
\end{gather*}
Hence, on the set 
$$\Delta_{n,m}:=\{\omega\in \Omega\, |\, \exists k\in [m, [n/r]-1]\, :  \, (Y_{kr}\cdots Y_{(k-1)r+1})(\omega)\in G_{1/n_0}\}
\, ,
$$
\begin{gather}\nonumber
\inf_{x,\, y\in S^+}\big(\log \langle y, A_n x\rangle 
-\log \|A_nx\|\big) \\
\label{inequality-coef}\ge 
-\log (n_0) +\min_{mr\le \ell\le n-1} \big(\log v( Y_{\ell +1}^t 
\cdots Y_n^t )-\log \|  Y_{\ell+1}^t 
\cdots Y_n^t\| \big) \, .
\end{gather}
Notice that all the above quantities are non positive and 
that $\min_{mr\le \ell\le n} \big(\log v( Y_{\ell +1}^t 
\cdots Y_n^t )-\log \|  Y_{\ell+1}^t 
\cdots Y_n^t\| \big)$ has the same law as 
$\min_{1\le \ell\le n-mr} \big(\log v( Y_{\ell}^t 
\cdots Y_{1}^t )-\log \|  Y_{\ell}^t 
\cdots Y_{1}^t\| \big)$.

\medskip

Notice also that $\P(\Delta_{n,m}^c)= \eta^{[n/r-m]}$ for some $0\le \eta <1$, for $n_0$ large enough.

Then, we conclude thanks to Lemma 2.1, using Proposition 
\ref{delta} and Lemma \ref{deviation} with $\tilde \mu$ and taking $m=[n/r]-C\log n$, with $C|\log \eta|>1/2$ 
(always true if $\eta =0$).
\hfill $\square$

\medskip

We shall now obtain the rate $O(1/\sqrt n)$ for the spectral radius and the coefficients under a much stronger condition.

\medskip

\begin{Theorem}\label{theo-strong-cond}
Let $p\in (2,3]$. Assume that $\mu$ is strictly contracting and admits a moment of order $p$. Assume that $s^2>0$. Assume that there exists $0< \delta\le 1$ such that 
$\mu^{*r}(G_\delta)=1$. Then the conclusion of Theorem \ref{berry-esseen}
holds with $\log\big(\inf_{x,y\in S^+}\langle x, A_n y)\rangle$ or $\log \kappa(A_n)$ instead of $\log \|A_n\|$.
\end{Theorem}
\noindent {\bf Proof.} By assumption, for every $n\ge r$ and $x\in S^+$, 
using that $\frac{A_n x}{\|A_nx\|}=
(Y_n\cdots Y_{n+1-r})\cdot (A_{n-r}x)$, we have, for every $x,y\in S^+$
$$
1\ge \frac{\langle y, A_nx \rangle}{\|A_n x\|}\ge 
\delta \quad \mbox{$\P$-a.s.}
$$
Then, the result follows from Theorem \ref{berry-esseen} and the fact that $\|A_n\|\ge \kappa(A_n)\ge
\inf_{x,y\in S^+}\langle x, A_n y\rangle$ .\hfill $\square$

\medskip

\medskip

We now give a condition that is equivalent to the  condition  $\mu^{*r}(G_\delta)>0$.  An equivalent condition, specific to the case of positive matrices (hence not valid in the general situation considered in Section \ref{cones}), has been obtained in \cite{GLX}, see their Lemma 2.1.

\medskip

For every $C>0$  and $0\le\gamma <1$, set
$$
G_{C,\gamma}:=\{g\in G\, :\, c(g)\le \gamma\mbox{ and }
\|g\|\le C v(g^t)\}\, .
$$
\begin{Lemma}
For every $0<\delta\le 1$, there exists $0\le \gamma<1$ and $C>0$ 
such that $G_\delta\subset G_{C,\gamma}$. Conversely, for 
every $0\le\gamma'<1$ and every $C'>0$ there exists 
$0 <\delta'\le 1$ such that $G_{C', \gamma'}\subset G_{\delta'}$. Hence, there exists $0<\delta \le 1$ such that $\mu(G_\delta)>0$ if and only if there exists $0\le \gamma<1$ and $C>0$ such that $\mu(G_{C,\gamma})>0$.
\end{Lemma}

\noindent {\bf Proof.} The proof relies on the following observations: for every $x\in S^+$, $\langle x, g 
e\rangle =\|g^tx\|$ and $\|g^tx\|/\|g\|\ge \langle x, g\cdot e\rangle/d \ge  \|g^tx\|/(d\|g\|).$

\medskip

Let $g\in G_\delta$, with $\delta>0$. By the previous computations, 
$\|g\|\le v(g^t)/\delta$. 

\smallskip

Let $x,y\in S^+$. Let us bound $d(g\cdot x, g\cdot y)$. 
For every $u\in S^+$, we have 
$$
\delta \langle u, g\cdot y\rangle \le \delta \le 
\langle u, g\cdot x\rangle\, .
$$
This implies that $m(g\cdot x,g\cdot y)\ge \delta$ (notice that then we must have $\delta\le 1$. Similarly, 
$m(g\cdot y,g\cdot x)\ge \delta$ and $d(g\cdot x, g\cdot y)
\le \frac{1-\delta^2}{1+\delta^2}=:\gamma<1$.
 So, $G_\delta \subset G_{1/\delta,\gamma}$. 
 
 \medskip
 Let $0\le \gamma<1$ and $C>0$. Let $g\in G_{C,\gamma}$. 
 Let $x,y\in S^+$. Notice that $m(g\cdot x,g\cdot y)\le 1$. 
 Hence, $\gamma \ge d(g\cdot x, g\cdot y)\ge 
 \frac{1-m(y,x) }{1+m(y,x)}$ and $m(y,x)\ge \frac{1-\gamma}{1+\gamma}$. We infer that $g\cdot y- \frac{1-\gamma}{1+\gamma}
 g \cdot x$ has non negative coordinates. Taking, 
 $x=e$, we see that for every $u\in S^+$, 
 $$
\langle u,g\cdot y\rangle \ge  \frac{1-\gamma}{1+\gamma}\, 
\langle u, g\cdot e\rangle \ge \frac{1-\gamma}{1+\gamma}\,
\|g^tx\|/(d\|g\|)\ge \frac{1-\gamma}{Cd(1+\gamma)}\, .
 $$
 \hfill $\square$





\medskip

\section{Regularity of the invariant measure}

\setcounter{equation}{0}

We prove here regularity properties of the invariant measure under various moment conditions.

\begin{Theorem}
Assume that $\tilde \mu$ is strictly contracting and  admits a moment of order $p\ge 1$. 
Then 
\begin{equation}\label{regularity-polynomial}
\int_{S^+}\sup_{y\in S^+}|\log \langle y, x\rangle|^p 
\, d\nu(x)<\infty\, .
\end{equation}
\end{Theorem}

\noindent {\bf Remark.} In the case of invertible matrices, Benoist and Quint \cite{BQ} proved that 
under a moment of order $p\ge 1$, $\sup_{y\in X}\int_X  
|\log \langle y, x\rangle|^{p-1} 
\, d\nu(x)<\infty\, .$

\noindent {\bf Proof.} It is standard that it suffices to prove that $\sum_{n\ge 1}n^{p-1}\nu (\sup_{y\in X}  
|\log \langle y, \cdot \rangle|\ge cn )<\infty $, 
for some $c>0$. Using that $\nu$ is $\mu$-invariant, it suffices to prove that 
$$
\sum_{n\ge 1} n^{p-1}\P\big(\sup_{x, \, y\in S^+}  \big|\log \langle y, A_n\cdot x\rangle \big|\ge cn\big)<\infty\, .
$$
Now, on $\Delta_{n,1}$,  by \eqref{inequality-coef},
\begin{equation}\label{delta_n_1}
\big|\log \langle y, A_n\cdot x\rangle \big| \le \log n_0+\max_{1\le k\le n} \big|\log v(Y_k^t\cdots Y_n^t)-\log \|Y_k^t\cdots Y_n^t\|\big|\, ,
\end{equation}
and we conclude thanks to Proposition \ref{delta} the fact that 
$\P(\Delta_{n,1}^c)\le \eta ^{[n/r-1]}$. \hfill $\square$

\medskip

\begin{Theorem}
Assume that $\tilde \mu$ is strictly contracting and  admits an exponential  moment of order $\gamma\in (0,1]$. 
Then, there exists $\delta>0$ such that
\begin{equation}\label{regularity-exponential}
\int_{S^+}\sup_{y\in S^+}{\rm e}^{\delta \big(-\log |\langle y, x\rangle|\big)^{\gamma} }
\, d\nu(x)<\infty\, .
\end{equation}
\end{Theorem}

\noindent {\bf Proof.} Proceeding as above, the theorem will be proved if we prove that there exist $\delta,\eta>0$ such that
\begin{equation}\label{regularity-expo}
\sum_{n\ge 1}{\rm e}^{\delta n^\gamma}\P(\sup_{x,y\in S^+}\big| \log \langle y, A_n \cdot x \rangle \big|\ge \eta n)<\infty\, .
\end{equation}
 We conclude thanks to \eqref{delta_n_1} and Lemma \ref{deviation}. \hfill $\square$

\medskip

\section{Deviation inequalities}

\setcounter{equation}{0}

\medskip



We now provide deviation estimates.

\begin{Proposition}\label{deviation-polynomial}
Assume that $\mu$ is strictly contracting and admits a moment of order $p\ge 1$. Let $\alpha\in (1/2, 1]$ such that $\alpha\ge 1/p$. For any $\varepsilon>0$, we have 
$$
\sum_{n\ge 1}n^{\alpha p-2}\sup_{x\in S^+}\P(\max_{1\le k\le n} |\sigma(Y_k,A_{k-1}\cdot x)-k\lambda_\mu|\ge 
n^\alpha \varepsilon)<\infty\, . 
$$
\end{Proposition}

\noindent {\bf Remark.} Using Proposition \ref{delta} and the fact that for $Z\in L^p$, $p\ge 1$,
$\sum_{n\ge 1} n^{p\alpha -1} \P(Z\ge n^\alpha \varepsilon)<\infty$, for any $\varepsilon>0$ and any $\alpha>0$, one can prove similar results for $\log \|A_n\|-n\lambda_\mu$, 
$\log \kappa(A_n)-n \lambda_\mu$, $\log v(A_n)-n\lambda_\mu$ or $\sup_{x\in S^+}|\log \|A_nx\|-n\lambda_\mu|$.  

\medskip

Proposition \ref{deviation-polynomial} is the version for positive matrices of  Theorem 4.1 of \cite{CDM-deviation}, stated for invertible matrices. The proof is exactly the same. Let us mention the key ingredients: 
The result concerns a cocycle for which, when $p\ge 2$, 
the function $\psi$ in \eqref{psi} is well defined and bounded and $\sup_{k\ge 1}\sup_{x\in S^+}\|\E((\sigma(Y_k, A_{k-1}\cdot x))^2
|{\mathcal F}_{k-1})\|_\infty<\infty$; and, when $1\le p<2$, one can control the coefficients 
$\delta_{1, \infty}(n)$.

\medskip

Concerning the matrix coefficients, the following result holds.

\medskip

\begin{Proposition}
Assume that $\mu$ is strictly contracting and that $\mu$ and $\tilde \mu$ admit a moment of order $p\ge 1$. Fro any $\varepsilon>0$,
Then
$$
\sum_{n\ge 1}n^{\alpha p-2}\P(\sup_{x, \, y\in S^+} |\log \langle y, A_n x\rangle -n\lambda_\mu|\ge 
n^\alpha \varepsilon)<\infty\, . 
$$
\end{Proposition}

\noindent {\bf Remark.} One cannot expect to have a maximum over $1\le k\le n$ inside the probability, since one may have 
$\P(\log \langle y, A_1 x\rangle=-\infty)>0 $, for some $x,\, y\in S^+$.

\noindent {\bf Proof.}  Using \eqref{inequality-coef} 
with $m=1$, we see that on $\Delta_{n,1}$
\begin{gather*}
\sup_{x,\, y\in S^+} |\log \langle y, A_n x\rangle 
-n\lambda_\mu|\\ \le \sup_{x\in S^+}|\log \|A_n x\|-n\lambda_\mu |+ \max_{1\le k\le n} \big|\log v(Y_k^t\cdots Y_n^t)-\log \|Y_k^t\cdots Y_n^t\|\big|.
\end{gather*}
To conclude, we apply then remark after Proposition \ref{deviation-polynomial} and the fact that the random variables $\max_{1\le k\le n} \big|\log v(Y_k^t\cdots Y_n^t)-\log \|Y_k^t\cdots Y_n^t\|\big|$ and
$\max_{1\le k\le n} \big|\log v(Y_k^t\cdots Y_1^t)-\log \|Y_k^t\cdots Y_1^t\|\big|$ have the same law, combined with Proposition \ref{delta} applied to $\tilde \mu$. \hfill $\square$

\section{Generalization to cones}\label{cones}

\setcounter{equation}{0}

In this section we show how to extend  the previous results to general cones. In the previous sections we studied products of positive matrices, that is products of matrices  
leaving  invariant the cone $(\R^+)^d$. In this section we consider more general cones.
This type of generalization was also investigated in \cite{BDGM}.

\medskip

There are many examples of closed solid cones as the ones considered below. For instance, the Lorentz (or ice-cream) 
cone: $\{(x_1, \ldots , x_n, z)\in \R^{n+1}\, :\, z\ge 0, \, 
x_1^2+\ldots +x_n^2\le z^2\}$. The linear operators (of matrices) leaving invariant the Lorentz cone have been studied in details by Loewy and Schneider \cite{LS}.

\medskip

Another example is the cone $K_S$ of positive semi-definite matrices of order $n$
viewed as a cone of the vector space of symmetric matrices of order $n$. 
Examples of operators leaving invariant $K_S$ are given by 
$M\mapsto A^tMA$ where $A$ is a matrix of size $n$ or
$M\mapsto {\rm tr} (MR_0)S_0$, with $R_0, \, S_0\in K_S$ and convex combinations of those.

\medskip

Let $d\ge 2$. We endow $V=\R^d$ with its usual inner product $\langle \cdot, \cdot\rangle$ and the 
associated norm $\|\cdot \|_2$. 


\medskip

Let $K$ be a closed proper convex cone with non empty interior
 of $\R^d$.   We recall that a cone of $\R^d$ is a set of $\R^d$ stable by multiplication by 
 non-negative real numbers and that it is proper if $K\cap(-K)=
 \{0\}$.  
 
We shall call such cones \emph{closed solid cones}. 

\medskip

Usually, the term \emph{solid cone}, refers only to a cone with non empty interior as in  \cite{LN}, 
 page 3. Hence, we add the convexity and the fact that $K\cap (-K)=\{0\}$.

 \medskip



 We associate with $K$ its dual cone $K^*:= 
 \{x^*\in V^*\, :\, \langle x^*,x\rangle\ge 0 \quad \forall x\in V\}$.

 \medskip
 
 By Lemma 1.2.4 of \cite{LN},  $K^*$ is also a closed solid  cone. Moreover, for every $x^*\in {\rm int }(K^*)$,  (the interior of 
 $K^*$)  $\langle x^*,x\rangle >0$ for every $x\in K\backslash \{0\}$ and $\Sigma_{x^*}:=\{x\in K\, :\, \langle x^*,x\rangle = 1
 \}$ is a compact convex set.
 
 \medskip

We define a partial order on 
 $V$ by setting for every $x,\, y\in V$, $x\preceq_K y$ if 
 $y-x\in K$.

 \medskip

 In the sequel we will need to work with a monotone norm for $K$, 
 that is a norm compatible with $\preceq_K$ in the sense of \eqref{monotone} below. 
 
\medskip

 Let us fix once and for all $x_0^*\in {\rm int}(K^*)$. Then, for every 
 $x\in V$, set 
 \begin{equation}\label{norm-monotone}
\|x\|_{x_0^*}=\sup _{x^*\in K^*\, :\, x^*\preceq_{K^*}x_0^*}\langle x^*,x\rangle\, .
 \end{equation}
 By Lemma \ref{lemme-norm}, $\|\cdot \|_{x_0^*}$ is a norm on $V$ and, using the definition of $K^*$, 
 \begin{equation}\label{monotone}
 \|x\|_{x_0*}\le \|y\|_{x_0^*} \qquad \forall 0\preceq_{K}
 x\preceq_K y\, .
 \end{equation}

 \medskip

 Notice also that 
 \begin{equation}\label{norm-simple}
 \|x\|_{x_0^*}=\langle x_0^*,x\rangle\qquad \forall x\in K\, .
 \end{equation}

 \medskip
 
 Recall that $(K^*)^*=K$. Hence fixing once and for all some $x_0\in {\rm int}(K)$, with $\langle x_0^*,x_0\rangle=1$, ,  one defines also 
 a monotone norm on $V^*$ by setting
 $$
\|x^*\|_{x_0}:= \sup_{x\preceq _K x_0}|\langle x^*, x\rangle |
\qquad \forall x^*\in V^*\, .
 $$ 
  Then, for every $x^*\in K^*$, $\|x^*\|_{x_0}=\langle x^*, x_0\rangle$.

\medskip

Set 
$$
S^{+}:=K\cap \{x\in V\, :\, \|x\|_{x_0^*}=1\}=\{x\in K\, :\, \langle x_0^*,x\rangle =1\}\, 
$$
and
$$
S^{++}:={\rm int}(K)\cap \{x\in V\, :\, \|x\|_{x_0^*}=1\}
=\{x\in {\rm int}(K)\, :\, \langle x_0^*,x\rangle =1\}\,  .
$$
 
 Notice that those definitions are consistent with 
 \eqref{S+-1} and \eqref{S+-2}, taking $x_0^*=(1,\ldots , 1)$.

 \medskip

 We shall now define an application $d$ on $(K\backslash\{0\})^2$ that will make $(S^+,d)$ a metric space.
 
 \medskip
 
 We first define an equivalence relation $\sim_K$ on $K$, by setting for every $x,\, y$, $x\sim_K y$ if there exists $0<\alpha \le \beta$ such that $\alpha x\preceq_K y\preceq \beta x$. The equivalence classes for $\sim_K$ are called \emph{parts} of $K$. 
By Lemma \ref{lemme-cone}, ${\rm int}(K)$  is a part of $K$.

 \medskip
 
 Given $x, \, y\in K\backslash\{0\}$, set 
 $$
m(x,y)=\sup\{\lambda \ge 0\, :\,  \lambda y\preceq_K x\}\, . 
 $$
 This definition is consistent with the definition of the function $m$ defined in Section 1 when $K=(\R^+)^d$.
 
 \smallskip

 Notice that if some $\lambda>0$ is such that $\lambda y
 \preceq_K x$ then $x-\lambda y\in K$, hence $x/\lambda-y \in K$. So $m(x,y)<+\infty$ since $K $ is closed and  $K\cap (-K)=\{0\}$.
 
 \medskip

 In particular, using again that $K$ is closed, $m(y,x)m(x,y)y\preceq_K m(y,x) x\preceq _K y$ so that $m(y,x)m(x,y)\le 1$.
 
 \medskip

 Then, we define for every $x,\, y\in K\backslash \{0\}$, 
 $$
d(x,y)=\varphi(m(x,y)m(y,x))\, ,
 $$
 where $\varphi$ is given by \eqref{phi}
 \medskip
 
  It follows from the definition of $\sim_K$ that 
  $x\sim_K y$ if and only if $m(x,y)m(y,x) =0$ if and only if 
  $d(x,y)=1$.

 Then, $d(x,y)= \tanh \big(\, (1/2) d_H(x,y)\, \big)$ where $d_H$ is introduced 
 page 26 of \cite{LN}. Actually, $d_H$ is only defined when 
 $x\sim_K y$ to avoid situations where $d_H(x,y)=+\infty$. 
 
 \medskip

\begin{Proposition}\label{complete-metric}
$(S^+,d)$ is a complete metric  space and $S^{++}$ is closed. Moreover, there exists $C_{x_0}>0$ such that 
\begin{equation}\label{norm-metric}
\|x-y\|_{x_0^*}\le C_{x_0^*} \frac{d(x,y)}{1-d(x,y)}\qquad \forall (x,y)\in \, S^{+}. 
\end{equation}
\end{Proposition} 

\noindent {\bf Remark.} When $x\sim_Ky$ the right-hand side 
of \eqref{norm-metric} is finite. Otherwise, $d(x,y)=1$ and the 
right-hand side of \eqref{norm-metric} has to be interpreted as $+\infty$.

\medskip

 \noindent {\bf Proof.} We first prove that $(S^+,d)$ is a  metric  space. Let $x,y,z\in S^+$ be 
 such that $x\sim_K y$ and $y\sim_K z$. By Proposition 2.1.1 of \cite{LN}, 
 $d_H(x,z)\le d_H(x,y)+d_H(y,z)$. Using that $u\mapsto \tanh (u/2)$ is subadditive, the inequality remains true with  $d$ in place of $d_H$.  If we do not have $x\sim_K y$ and $y\sim_K z$, then $m(x,y)m(y,x)=0$ or $m(y,z)m(z,y)=0$, hence $d(x,y)=1$ or $d(y,z)=1$ so that the triangle inequality is still satisfied. 
 
 The fact that $d$ is a distance on $S^+$ then follows from (other statements of) Proposition 2.1.1 of \cite{LN}. The fact that 
 $(S^{+},d)$ is complete follows 
 from Lemma 2.5.4 of \cite{LN}. Indeed, if $(x_n)_{n\in \N}\subset S^+$ is a Cauchy sequence for $d$, then $d(x_p,x_q)<1$, say for 
 $q,\, p\ge N$, so that $(x_n)_{n\ge N}$ is included in a part $P$ of $K$. But, by Lemma 2.5.4 of \cite{LN}, $S^+\cap P$ is complete for $d$. 
 
 \medskip
 
 Let us explain why $S^{++}$ is closed. Using similar arguments as above we see that it is enough to prove that ${\rm int}(K)$ is a part of $K$, but this follows from Lemma 
 \ref{lemme-cone}.

 \medskip

 \eqref{norm-metric} follows from 
 (2.21) page 47 of \cite{LN}, 
 using the relation between $d_H$ and $d$. \hfill $\square$

 \medskip

 \medskip
 
 We shall now define the analogue of the positive matrices.
 
 \medskip
 
 \medskip

 Let 
 $$G:=\{g\in M_d(\R)\, :\, g(K\backslash\{0\})\subset K\backslash\{0\},\, g({\rm int}(K))\subset {\rm int}(K)
 \}\, .
 $$ 
 
 \medskip

It follows from Lemma \ref{lemme-allowable} that 
$$G:=\{g\in M_d(\R)\, :\, g^t(K^*\backslash\{0\})\subset K^*\backslash\{0\},\, g^t({\rm int}(K^*))\subset {\rm int}(K^*)
 \}\, .
 $$ 
In particular, $g\in G$ is  allowable in the sense of \cite{BDGM} (see $a)$ page 1527). Hence, the allowability condition in \cite{BDGM} is redundant.
 
 \medskip

 We endow $M_d(\R)$ with the norm: $\|g\|_{x_0^*}:=\sup_{x\in K,  \, \|x\|_{x_0^*}=1}\|gx\|_{x_0^*}$. The fact that this is indeed a norm follows from the fact that $K$ has non empty interior (i.e.  $K-K=V$). Notice that for $g\in G$, 
 $$
 \|g\|_{x_0^*}=\sup_{x\in K,  \, \langle x_0^*,x\rangle=1}
 \langle x_0^*, gx\rangle \, .
 $$

 \medskip
 
 Define also 
 $$G^+:=\{g\in G\, :\, g(K\backslash\{0\})\subset 
 {\rm int}(K)\}\, .
 $$
 
 By Lemma 10.1, 
 $$G^+:=\{g\in G\, :\, g^t(K^*\backslash\{0\})\subset  {\rm int}(K^*)\}\, .
 $$

 \medskip

 \medskip
 
 Define for every $g\in G$
 $$
v_{x_0^*}(g)=\inf_{x\in K,  \, \|x\|_{x_0^*}=1}\|gx\|_{x_0^*}\, ,  
 $$
 
 Notice that for $g\in G$, $v(g)=\inf_{x\in K,  \, \langle x_0^*,x\rangle=1}
 \langle x_0^*, gx\rangle $.
 
 \medskip
 
We then define $N_{x_0^*}(g):=\max(\|g\|_{x_0^*}, 1/v_{x_0^*}(g))$ and $L_{x_0^*}(g):= \frac{\|g\|_{x_0^*}}{v_{x_0^*}(g)}$. 
 
 \medskip
 
 The semi-group $G$ is acting on $S^+$ as follows.
 $$
g\cdot x =\frac{gx}{\|gx\|_{x_0^*}}=
\frac{gx}{\langle x_0^*, gx\rangle} \qquad \forall (g,x)\in G\times 
S^+\, .
$$

We then define a cocyle by setting $\sigma(g,x)=
\log (\|gx\|_{x_0^*})$ for every $(g,x)\in G\times S^+$.

\medskip

For every $g\in G$ set 
$$
c(g):= \sup_{x,\, y\in K\backslash\{0\}}d(g x, g y)\, .
$$

 \begin{Proposition}\label{prop-hennion-general}
For every $(g,g',x,y)\in G^2\times(S^+)^2$ we have
\begin{itemize}
\item [$(i)$] $|\sigma(g,x)|\le \log N(g)$;
\item [$(ii)$] $|\sigma(g,x)-\sigma(g,y)|\le 2C_{x_0^*}L(g)d(x,y)$ if $d(x,y)\le 1/2$;
\item [$(iii)$] $|\sigma(g,x)-\sigma(g,y)|\le 2\ln \big(1/
(1-d(x,y))\big)$ ;
\item [$(iv)$] $c(gg')\le c(g)c(g')$;
\item [$(v)$] $c(g)\le 1$ and $c(g)< 1$ 
iff $g\in G^+$;
\item [$(vi)$] $d(g\cdot x,g\cdot y)\le c(g) d(x,y)$.
\end{itemize}
\end{Proposition}

\noindent {\bf Remark.} The constant $C>0$ appearing in item $(ii)$ is the same as in \eqref{norm-metric}.

\noindent {\bf Proof.} $(i)$ is obvious. $(ii)$ may be proved exactly as item $(i)$ of Lemma 5.3 of \cite{Hennion}, using \eqref{norm-metric}. 

\medskip

Let us prove $(iii)$. Let $x,y\in S^+$. Assume that $x\sim_K y$, 
since otherwise the right-hand side in item $(iii)$ equals $+\infty$ and the inequality is clear. We have 
$m(x,y)y\preceq_K x $ and $m(y,x)x\preceq_K y $. Since $g\in G$, 
$m(x,y)gy\preceq_K gx $ and $m(y,x)gx\preceq_K gy $. Using that 
$\|\cdot \|_{x_0^* }$ is monotone we infer that $m(x,y)\|gy\|_{x_0^*}
\le \|gx\|_{x_0^*}$ and $m(y,x)\|gx\|_{x_0^*}
\le \|y\|_{gx_0^*}$. Hence 
$$
m(x,y)\le \frac{\|gx\|_{x_0^* }}{\|y\|_{x_0^*}}\le 1/m(y,x)\, .
$$
 Then, the proof may be finished as the proof of item $(ii)$ of Lemma 5.3 of \cite{Hennion}.
 
 \medskip
 
 The proof of $(iv)$ may be done exactly as in \cite{Hennion}. For the proof of $(v)$ we need to check some of the arguments. 
 
 \medskip
 
 Let $g\in G^{+}$. Then, $gS^+$ is a compact set (for $\|\cdot \|_{x_0^*}$) of ${\rm int}(K)$. Let us prove that is also compact for $d$.  Let $(x_n)_{n\in \N}\subset S^+$. Taking a subsequence if necessary, we may assume that there exists $y\in {\rm int}(K)$ such that $(gx_n)_{n\in \N}$ converges for $\|\cdot \|_{x_0}$ to 
 $y$. Since $y\in {\rm int}(K)$, by Lemma 2.5.5 of \cite{LN}, 
 $(x_n)_{n\in \N}$ converges to $y$ for $d_H$, hence for $d$. 
 
 The rest of the proof is as in \cite{Hennion}.
 
 Item $(vi)$ is just Birkhoff's inequality, see for instance 
 page 31 of \cite{LN}. \hfill $\square$
 
 \medskip
 
 We shall now consider the analogous statements as those given in Lemma \ref{properties}. Only item $(ii)$ requires a proof. 
 
 \begin{Lemma}
 There exists $C>0$ such that for every $g\in G$, 
 $$
\|gx_0\|_{x_0^*}\le \|g\|_{x_0^*}\le C \|gx_0\|_{x_0^*}\, .  
 $$
 \end{Lemma}
 
 \noindent  {\bf Proof.} Since $\langle x_0^*,x_0\rangle=1$, 
 $\|gx_0\|_{x_0^*}\le \|g\|_{x_0^*}$. 
 
 \medskip
 Let $x\in K$ be such that 
 $\langle x_0^*,x\rangle =1$. Let $g\in G$
 
 \medskip
 
 Using Lemma \ref{lemme-cone} with the cone $K^*$  there exists $\varepsilon >0$  such that 
 $g^t x_0^*\preceq_{K^*} 
 \frac{\|g^tx_0^*\|_{x_0}}\varepsilon x_0^*$. Hence, using 
 that $gx\in K$ and Lemma \ref{lemme-LN},
 \begin{gather*}
 \|gx\|_{x_0^*}= \langle x_0^*,gx\rangle =\langle g^tx_0^*
 ,x\rangle \le \frac{\|g^tx_0^*\|_{x_0}}\varepsilon\langle x_0^*,x\rangle
 \\
 =\frac{\langle g^tx_0^*,x_0\rangle}\varepsilon=\frac{\langle x_0^*, gx_0\rangle}\varepsilon
 =\frac{\|gx_0\|_{x_0^*}}\varepsilon\, .
 \end{gather*}
 \hfill $\square$
 
 All the results of the previous sections 
 hold true for a cocycle satisfying all the properties listed in Proposition \ref{prop-hennion} and Lemma \ref{properties}.

 \section{Technical results}
 
 \setcounter{equation}{0}
 
 The next lemma is just  Lemma 1.2.4 of 
 \cite{LN}. 
 
\begin{Lemma}\label{lemme-LN}
 Let $K$ be a closed solid cone. Then 
 $${\rm int} (K^*)=\{x^*\in V^*\, :\, \langle x^*, x\rangle >0\, ,\, \forall 
 x\in K\backslash\{0\}\}\, . $$
 In particular, 
 \end{Lemma} 
 
The next lemma follows from the proof Lemma 1.2.4 of 
 \cite{LN}.  We recall the short argument.
 
 \begin{Lemma}\label{lemme-cone}
Let $\|\cdot \|$ be a norm on $V=\R^d$. Let $K $ be a closed solid cone. Then, 
 for every $x\in {\rm int} (K)$, there exists $\varepsilon>0$, such that for every $y\in K\cap \bar B(0,1)$, where 
$\bar B(0,1)$ is the closure of the  unit ball $B(0,1)$, we have $y\preceq \frac1\varepsilon x$. Then $\|y\|\le \frac1\varepsilon$. In particular, ${\rm int}(K)$ is a part of $K$.
 \end{Lemma}
 
 \noindent {\bf Proof.} Let $x\in {\rm int} (K)$. There exists 
 $\varepsilon>0$ such that $\bar B(x,\varepsilon)\subset 
 {\rm int}(K)$. Let $y\in \bar B_{\|\cdot \|}(0,1)$. 
 Then, $x-\varepsilon y\in K$, which means precisely that 
 $y\preceq \frac1\varepsilon x$. In particular, if 
 $x,\, y\in {\rm int}(K)$, $x\sim_K y$. 
 
 \medskip
 
 It remains to prove that for every $(x,y)\in {\rm int}(K=\times K$, $x\sim_K y\Rightarrow y\in {\rm int}(K)$.
 
 \medskip
 
  Hence, let $x\in {\rm int}(K)$. There exists $\varepsilon >0$ such $B(x,\varepsilon)\subset K$.
 
 \smallskip
 
  Let $y\in K$ be such that $y\sim_K x$. There exists 
  $\alpha>0$ such that $x\preceq_K \alpha y$. So $\alpha y-x\in K$  and 
 $$
\alpha y=  x+ \alpha y-x\in \cup_{z\in K}(z+B(x,\varepsilon))\, ,
 $$
 which is an open subset of $K$.\hfill $\square$

 \begin{Lemma}\label{lemme-allowable}
 Let $g\in M_d(\R)$ and let $K$ be a closed solid cone of $E$. 
 \begin{itemize}
\item [$(i)$]  $g (K\backslash \{0\})\subset K\backslash \{0\}$ if and only if $g^t ({\rm int}(K^*))\subset {\rm int}(K^*)$;
\item [$(ii)$] $g ({\rm int}(K))\subset {\rm int}(K)$ if and only if $g^t (K^*\backslash \{0\})\subset K^*\backslash \{0\}$.
\end{itemize}
 \end{Lemma}
 
 \noindent {\bf Proof.} Assume that $g (K\backslash \{0\})\subset K\backslash \{0\}$. Let $x^*\in {\rm int}(K^*)$ and $x\in 
 K\backslash \{0\}$. We have
 \begin{gather*}
 \langle g^tx^*, x\rangle =\langle x^*, gx\rangle>0\, ,
 \end{gather*} 
by Lemma   \ref{lemme-LN}. Using Lemma \ref{lemme-LN} again, 
we see that $g^tx^* \in {\rm int}(K^*)$. 

\medskip

Assume that $g^t ({\rm int}(K^*))\subset {\rm int}(K^*)$. 
Let $x\in K\backslash \{0\}$ and $x^*\in {\rm int} (K^*)$. We have 
$$
\langle x^*, gx\rangle =\langle g^t x^*, x\rangle > 0\, .
$$
Hence $gx\in K^{**}=K$ (see Exercise 2.31 of \cite{BV}) and $gx\neq 0$, which proves item $(i)$. 

\medskip

Item $(ii)$ is just item $(i)$ for $K^*$ using that $K^{**}=K$. 
\hfill $\square$

\begin{Lemma}\label{lemme-norm}
$\|\cdot \|_{x_0^*}$ defined by \eqref{norm-monotone} is a norm for every $x_0^*\in {\rm int}(K^*)$.
\end{Lemma}

\noindent {\bf Proof.}  By Lemma 1.2.5 of \cite{LN}, the set $\{x^*\in K\, :\, x^*\preceq_{K^*}x_0^*\}$ is bounded, hence $\|\cdot \|_{x_0^*}$ is finite on $V$. The fact that $\|\cdot \|_{x_0^*}$ satisfies the triangular inequality and is positively homogeneous are obvious. 
 
 \medskip

 Assume that $x\in E$, is such that $\|x\|_{x_0^*}=0$. By Lemma 
 \ref{lemme-cone} applied to $K^*$ (with $x=x_0^*$),  
 for every $x^* \in K^*$, $\langle x^*,x\rangle=0$. Since $K^*$ has non empty interior, $K^*-K^*=V^*$ and $x=0$.  \hfill $\square$

\end{document}